 \documentclass{jfa}
\begin{document}

%


\authorrunninghead{Andrea Posilicano}

\titlerunninghead{A Krein-like formula}





\newcommand{\der}{\frac{d{\ }}{dz}\,}
\newcommand{\calH}{\mathcal H}
\newcommand{\X}{\mathcal X}
\newcommand{\Y}{\mathcal Y}
\newcommand{\Z}{\mathcal Z}
\newcommand{\Q}{\mathcal Q}
\newcommand{\R}{\mathcal R}
\newcommand{\G}{\mathcal G}
\newcommand{\K}{\mathcal K}
\newcommand{\E}{\mathcal E}
\newcommand{\F}{\mathcal F}
\newcommand{\M}{\mathcal M}
\newcommand{\RE}{\mathbf R}
\newcommand{\C}{{\mathbf C}}
\newcommand{\GB}{\breve G(z)}
\newcommand{\LD}{L^2(\RE^3)}
\newcommand{\sob}[2]{H^{#1}(\RE^{#2})}
\newcommand{\sobc}[1]{H^{#1}(C)}
\newcommand{\ld}[1]{L^2(\RE^#1)}
\newcommand{\HD}{H^{2}(\RE^3)}
\newcommand{\p}{\par\noindent}
\newcommand{\vp}{\varphi}
\newcommand{\square}{\vbox{\hrule\hbox{\vrule\vbox to 7 pt {\vfill\hbox to
                     7 pt {\hfill\hfill}\vfill}\vrule}\hrule}}
\newcommand{\re}{{\mathrm {Re}}}
\newcommand{\im}{{\mathrm {Im}}}
\newcommand{\uno}{I}
\newcommand{\wtilde}{\widetilde}

\title{A KREIN-LIKE FORMULA FOR SINGULAR PERTURBATIONS OF 
       SELF-ADJOINT OPERATORS AND APPLICATIONS}

\author{Andrea Posilicano}

\affil{Dipartimento di Scienze, Universit\`a dell'Insubria, I-22100 Como, Italy}

\email{posilicano@mat.unimi.it}

\abstract{Given a self-adjoint
operator $A:D(A)\subseteq\calH\to\calH$ and a continuous
linear operator $\tau:D(A)\to\X$ with Range$\,\tau'\cap\calH'
=\left\{0\right\}$, $\X$ a Banach space, we
explicitly construct a family $A^\tau_\Theta$ of self-adjoint operators such that any
$A^\tau_\Theta$ coincides with the original $A$ on the kernel of
$\tau$. Such a family is obtained 
by giving a Kre\u\i n-like formula where the role of the deficiency
spaces is played by the dual pair $(\X,\X')$; the parameter 
$\Theta$ belongs to the space of symmetric operators from $\X'$ to
$\X$. When $\X=\C$ one
recovers the ``$\,\calH_{-2}$ -construction'' of Kiselev and Simon and
so, to some extent, our results can be regarded as an extension
of it to the infinite rank case. Considering the situation in which
$\calH=L^2(\RE^n)$ and $\tau$ is the trace (restriction)
operator along some null subset, we give various applications to singular
perturbations of non necessarily elliptic pseudo-differential operators,
thus unifying and extending previously known results.}


\begin{article}

\section{INTRODUCTION}
Let $A:D(A)\subseteq\calH\to\calH$ be a self-adjoint operator on the Hilbert space
$\calH$ and suppose that there exists a linear dense set $N\subset D(A)$ which is 
closed with respect to the graph
norm on $D(A)$. If we denote by $A_N$ the restriction of $A$ to $N$, then $A_N$
is a closed, densely defined, symmetric operator. Since $N\neq D(A)$, $A$ is a non-trivial extension of $A_N$ and so, by the von Neumann theory on
self-adjoint extensions of closed symmetric operators (see \cite{[N]},
[17, \S XII.4], [35, \S X.1]), 
we know that the
deficiency indices $n_\pm$, defined as the dimensions of
$K_\pm:=$Kernel$\,A_N^*\pm i$, are equal and strictly positive. The family of 
self-adjoint extensions of $A_N$ is
then parametrized by the unitary maps from $K_+$ onto $K_-$. When $A$
is strictly positive, a deeper and more explicit construction of the
(positive if 
dim$\,K=+\infty$, $K:=$Kernel$\,A_N^*$) self-adjoint extensions of $A_N$
is given by the Birman-Kre\u\i n-Vishik theory (see \cite{[K3]}, \cite{[V]}, \cite{[B]},
\cite{[AS]}). In this case the family of (positive) extensions is parametrized by
the (positive) quadratic forms on $K$. \par
Any self-adjoint extension $\wtilde A_N\not= A$ can then be interpreted as a 
singular perturbation of $A$
since the two operators 
differ only on $\calH\backslash N$, the set $\calH\backslash N$ being ``thin'' since its complement is a 
linear dense subset of $\calH$.\par
In the case $n_{\pm}=1$, Kre\u\i n obtained, in 1943 (see \cite{[K1]}), a quite
explicit formula relating the resolvents of any two self-adjoint
extensions of a given symmetric operator. Such a formula was then extended,
by Kre\u\i n himself in 1946 (see \cite{[K2]}), to the case $n_\pm=m<+
\infty$. In our setting it states the following: 
for any $z\in\rho(A)\cap\rho(\wtilde A_N)$ one has 
$$(-\wtilde A_N+z)^{-1}=(-A+z)^{-1}
+\sum_{j,k=1}^m\Gamma(z)^{-1}_{jk}\,\varphi_j(z)\otimes\varphi_k(z^*)\, ,
$$
where 
$$
\varphi_k(z):=\varphi_k-(i-z)(-A+z)^{-1}\varphi_k\ ,
$$
$\left\{\varphi_k\right\}_1^m$ being the set of linear independent solutions of 
$$
A^*_N\varphi=i\,\varphi\,,\qquad\varphi\in D(A^*_N)\ ,
$$
and where the invertible matrix $\Gamma(z)$ satisfies
($\langle\cdot,\cdot\rangle$ denoting the scalar product on $\calH$)
$$
\Gamma(z)_{jk}-\Gamma(w)_{jk}
=(z-w)\langle\varphi_j(z^*),\varphi_k(z)\rangle\ .
$$
By such a formula, since $N$ is dense, one can then readily define $\wtilde A_N$ as
$$
D(\wtilde A_N):=\{\phi\in\calH\ :\
\phi=\phi_z+\sum_{j,k=1}^{m}\Gamma(z)_{jk}^{-1}\langle\varphi_k(z),
\phi_z\rangle\,\varphi_j(z)\,,\ \phi_z\in D(A)\,\}
$$
$$
(-\wtilde A_N+z)\phi:=(-A+z)\phi_z\ .
$$
Kre\u\i n's original papers were written in russian, but his results were
popularized in some excellent monographs (see
e.g. [1, chap. VII]). Instead, the analogous formula for the case
$n_\pm=+\infty$, which was obtained by Saakjan in 1965 (see \cite{[S]}), is much less
known, since the work is not available in english (see however \cite{[GMT]} and
references therein). Due probably to this fact, the Kre\u\i n formula for 
$n_\pm=+\infty$ (similar considerations also apply to the
Birman-Kre\u\i n-Vishik theory) was rarely used in concrete applications: we are 
mainly referring to the much studied case of singular perturbations 
of the Laplacian supported by null sets (see e.g. \cite{[AGHH]},
\cite{[AFHKL]}, 
\cite{[Br]} and references 
therein). Indeed in situations of this kind other approaches are used:
extensions are mainly obtained
either as resolvent limits of less singular perturbations or by other
constructions often resembling variations of either the Kre\u\i n formula
or the Birman-Kre\u\i n-Vishik theory. Usually such approaches rely on the
elliptic nature of the Laplacian and are not applicable to the study of
singular perturbations of hyperbolic operators (this was the original 
motivation of our work).\par 
Here we show how, 
when the (not necessarily dense) set $N$ is the kernel of a
continuous linear map 
$\tau: D(A)\to \X$  such that
Range$\,\tau'\cap\calH'=
\left\{0\right\}$, $\X$ a Banach space, one can prove, by almost
straightforward 
arguments, a Kre\u\i n-like formula for a family $A^\tau_\Theta$,
$\Theta$ a symmetric operator from $\X'$ to $\X$, of self-adjoint extensions of $A_N$,
where the role of $K_\pm$ is played by the dual pair $(\X,\X')$ (our
construction could be given for $\X$ a locally convex space, but we
will not strive here for the maximum of generality). \par In contrast
to other approaches (see e.g. \cite{[S]}, \cite{[GMT]},
\cite{[DM1]}, \cite{[DM2]} 
and references therein) the
formula given here turns out to be relatively simple being expressed
directly in terms of the map $\tau$; moreover we do
not need to compute $A_N^*$. In more detail (see theorem 2.1) one
obtains, under a hypothesis which we prove to be satisfied under relatively
weak conditions (see proposition 2.1),
$$
(-A^\tau_\Theta+z)^{-1}=(-A+z)^{-1}+G(z)\cdot(\Theta+\Gamma(z))^{-1}\cdot\breve
G(z)\ ,
$$ 
where
$$
\breve G(z):=\tau\cdot (-A+z)^{-1}\,,\qquad
G(z):=C_\calH^{-1}\cdot\breve G(z^*)'
$$
($C_\calH$ being the canonical isomorphism of $\calH$ onto $\calH'$)
and the conjugate linear operator $\Gamma(z):D\subseteq\X'\to
\X$ 
satisfies the equation
$$
\forall\, \ell\in D,\qquad
\der \Gamma(z)\ell=\breve G(z)\cdot
G(z)\ell
$$
which (see lemma 2.2) we show to have an explicit (in terms of $\tau$
itself) bounded operator solution. Such a solution plays a fundametal role
in finding (see lemmata 2.3 and 2.4) other nicer (even if unbounded) 
solutions which we then use
in (some of) the examples.\par
In \S 3, after showing (example 3.1)
how our construction, in the case $\X=\C$, reproduces the
``$\,\calH_{-2}$ -construction'' given in \cite{[KS]} and how, in the case $A$ is
strictly positive, it gives a variation on the Birman-Kre\u\i n-Vishik
theory which comprises the results in \cite{[KKO]} (example 3.2), we use the above 
Kre\u\i n-like formula to study singular perturbations of
non necessarily elliptic
pseudo-differential operators, thus unifying and extending 
previously known results. More precisely we give the following
examples:\p
\begin{itemize}
\item Finitely many point interaction in three dimensions
(example 3.3);
\item Infinitely many point interaction in three dimensions (example 3.4);
\item Singular perturbations of the Laplacian in three and four dimensions
supported by regular curves (example 3.5);
\item Singular perturbations, supported by null sets with Hausdorff
codimension less than $2s$, of translation
invariant pseudo-differential
operators with domain $H^{s}(\RE^n)$ (example 3.6);
\item Singular perturbations of the d'Alembertian in four dimensions
supported by time-like straight lines (example 3.7). In order to limit the
lenght of the
paper we content ourselves with discussing here only the case of
a straight line. A complete study of the case of a generic time-like
curve will be the subject of a separate paper. We belive that the
detailed study of such a kind of operators will lead to a rigorous
framework for the classical and
quantum electrodynamics of point particles in the spirit of
the results obtained, for the linearized (or dipole) case, in \cite{[NP1]}-\cite{[NP3]}
and \cite{[BNP]};
\item Singular perturbations, supported by null sets, of translation
invariant pseudo-differential
operators with domain the Malgrange spaces $H_\vp(\RE^n)$ (example 3.8).
\end{itemize}

\subsection*{Definitions and notations}

\begin{itemize}
\item Given a Banach space $\X$ we denote by $\X'$ its strong dual;  
\item $L(\X,\Y)$, resp. $\wtilde L(\X,\Y)$, denotes the
space of linear, resp. conjugate linear, operators 
from the Banach space $\X$ to the Banach space $\Y$.
\item $B(\X,\Y)$, resp. $\wtilde B(\X,\Y)$, denotes the space of 
bounded, everywhere defined, linear, resp. conjugate linear, operators 
on the Banach space $\X$ to the Banach space $\Y$. It is a Banach space with the
norm $\|A\|_{\X,\Y}:=\sup\{\|Ax\|_\Y,\, \|x\|_\X=1\}$. 
\item The closed linear operator operator $A'$ and the 
conjugate linear closed operator $\wtilde A'$ are the adjoints of the densely 
defined linear operator $A$ and of the densely defined conjugate linear 
operator $\wtilde A$ respectively, i.e. \p
\centerline{\mbox {$\forall\,x\in D(A)\subseteq\X,\quad\forall
\ell\in D(A')\subseteq
\Y',\qquad (A'\ell)(x)=\ell(Ax)$,}}\p
\centerline{\mbox {$\forall\,x\in D(\wtilde
A)\subseteq\X,\quad \forall\ell\in D(\wtilde
A')\subseteq \Y',\qquad(\wtilde A'\ell)(x)=(\,\ell(\wtilde Ax)\,)^*$,}}\ 
where $^*$ denotes complex conjugation. 
\item $\wtilde S(\X',\X)$ denotes
the space of conjugate linear operators $A$ such that \p
\centerline{\mbox{$\forall\,\ell_1,\ell_2\in D(A),\qquad
\ell_1(A\ell_2)=(\,\ell_2(A\ell_1)\,)^*$. }}
\item For any $A\in\wtilde S(\X',\X)$ we define\p \centerline{\mbox{$\gamma(A):=
\inf\left\{\,\ell(A \ell),\ \|\ell\|_{\X'}=1,\ \ell\in D(A)\,\right\}$.}}\p
\item $J_\X\in B(\X,\X'')$ indicates the injective map (an
isomorphism when $\X$ is reflexive) defined by $(J_\X x)(\ell):=\ell(x)$.
\item If $\calH$ is a complex Hilbert space with scalar product 
(conjugate linear w.r.t. the first variable)
$\langle\cdot,\cdot\rangle$, then $C_\calH\in\wtilde B(\calH,\calH')$ denotes the isomorphism
defined by $(C_\calH y)(x):=\langle y,x\rangle$. The Hilbert
adjoint of the densely defined linear operator $A$ is then given by $A^*=
C_\calH^{-1}\cdot A'\cdot C_\calH$. 
\item $\F$ and $*$ denote Fourier transform and convolution
respectively.
\item $H^s(\RE^n)$, $s\in\RE$, is the usual scale of
Sobolev-Hilbert spaces, i.e. $H^s(\RE^n)$ is the space of tempered
distributions with a Fourier transform which is square integrable
w.r.t. the measure with density $(1+|x|^2)^s$.
\item $c$ denotes a generic strictly positive constant
which can change from line to line.
\end{itemize}

\section{A KREIN-LIKE FORMULA}
Let $$A:D(A)\subseteq\calH\to\calH$$ be a self-adjoint operator on the complex Hilbert
space $\calH$. $D(A)$
inherits a Hilbert space structure by introducing the usual scalar
product leading to the graph norm
$\|\phi\|^2_A:=\langle\phi,\phi\rangle+\langle
A\phi,A\phi\rangle$. Denoting the 
resolvent set of $A$ by $\rho(A)$ we define, for any $z\in\rho(A)$,
$$
R(z):=(-A+z)^{-1}:\calH\to D(A)\,,\qquad R(z)\in B(\calH,D(A))\ .
$$
We consider now a linear operator
$$
\tau:D(A)\to\X \,,\qquad \tau\in B(D(A),\X)\ ,
$$
where $\X$ is a complex Banach space. By means of $A$ and $\tau$ we
can define, for any $z\in\rho(A)$, the following operators:
$$
\breve G(z):=\tau\cdot R(z) :\calH\to\X\,,\qquad\breve G(z)\in B(\calH,\X)\ ,
$$
$$
G(z):=C_\calH^{-1}\cdot\breve G(z^*)':\X'\to\calH\,,\qquad
G(z)\in\wtilde B(\X',\calH)\ .
$$
\begin{remark}
As $R(z)$ is surjective, $R(z)'$ is injective.
If $\tau$ has dense range then $\tau'$ is
injective. Therefore, when $\tau$ has dense range, $\breve G(z)$ has
dense range and $G(z)$ is injective. This implies that the only
$\Lambda\in B(\X,\X')$ 
which solves the operator equation $G(z)\cdot\Lambda\cdot \breve G(z)=0$
is the zero operator.
\end{remark} 
\begin{lemma} For any $w$ and $z$ in $\rho(A)$ one has
\begin{eqnarray*}
(z-w)\,\breve G(w)\cdot R(z)&=&\breve G(w)-\breve G(z)\\
(z-w)\, R(w)\cdot G(z)&=&G(w)-G(z)\ .
\end{eqnarray*}
\end{lemma}
\begin{proof} By the first resolvent identity one has
$$(z-w)\,R(w)\cdot R(z)=R(w)-R(z)\ .
$$
Therefore
$$
(z-w)\,\breve G(w)\cdot R(z)=
(z-w)\,\tau\cdot R(w)\cdot R(z)
=\breve G(w)-\breve G(z)
$$
and, by duality (here $R(z)$ is considered as an element of $B(\calH,\calH)$),
\begin{eqnarray*}
G(w)-G(z)&=&C_\calH^{-1}\cdot\left(\breve G(w^*)-\breve G(z^*)\right)'\\
&=&C_\calH^{-1}\cdot\left((z^*-w^*)\,\breve G(z^*)\cdot R(w^*)\right)'\\
&=&C_\calH^{-1}\cdot\left((z^*-w^*)\,C_\calH\cdot R(w)\cdot C_\calH^{-1}\cdot
\breve G(z^*)'\right)\\
&=&(z-w)\, R(w)\cdot G(z)
\ .
\end{eqnarray*}
This ends the proof.
\end{proof}
\begin{remark} The second relation in the lemma above shows that
$$
\forall\,w,z\in\rho(A),\qquad \mbox{\rm
Range}\,(\,G(w)-G(z)\,)\subseteq 
D(A)\ . 
$$
\end{remark}
We want now to define a new self-adjoint operator which, when
restricted to the kernel of $\tau$, coincides with the original
$A$. Since, in the case of a bounded perturbation $V$, 
for any $z$ such that $\|V\cdot R(z)\|_{\calH,\calH}<1$ one has 
$$
(-(A+V)+z)^{-1}=R(z)+ R(z)\cdot(\,\uno-V\cdot
R(z)\,)^{-1}\cdot V\cdot R(z)\ ,
$$  
we are lead to write the presumed resolvent as 
$$
R^\tau(z)=R(z)+B(z)\cdot\tau\cdot R(z)\equiv R(z)+B(z)\cdot\breve G(z)\,, 
$$
where $B(z)\in B(\X,\calH)$ has to be determined. \par
Self-adjointness requires $R^\tau(z)^*=R^\tau(z^*)$ or, equivalently, 
\begin{equation}
G(z)\cdot B(z^*)'\cdot C_\calH=B(z)\cdot\breve G(z)\ .
\end{equation}
Therefore if we put $B(z)=G(z)\cdot \Lambda(z)$,
$\Lambda(z)\in\wtilde B(\X,\X')$, then one can check that (1) is implied
by 
(by remark 2.1, when $\tau$ has dense
range, is equivalent to)
\begin{equation}
\Lambda(z)'\cdot J_\X=\Lambda(z^*)\ .
\end{equation}
We now impose the resolvent identity 
\begin{equation}
(z-w)\,R^\tau(w)R^\tau(z)=R^\tau(w)-R^\tau(z)\ .
\end{equation}
Since (we make use of lemma 2.1)
\begin{eqnarray*}
&{\quad}&(z-w)\,R^\tau(w)\cdot R^\tau(z)=(z-w)\,\left(R(w)\cdot R(z)\right.\\
&{\quad}&+ R(w)\cdot G(z)\cdot\Lambda(z)\cdot\breve G(z)
+G(w)\cdot\Lambda(w)\cdot\breve G(w)\cdot R(z)\\
&{\quad}&+\left.G(w)\cdot\Lambda(w)\cdot\breve G(w)\cdot
G(z)\cdot\Lambda(z)\cdot\breve G(z)\right)\\
&=& R(w)- R(z)
+G(w)\cdot\Lambda(z)\cdot\breve G(z)
-G(z)\cdot\Lambda(z)\cdot\breve G(z)\\
&{\quad}&+G(w)\cdot\Lambda(w)\cdot\breve G(w)
-G(w)\cdot\Lambda(w)\cdot\breve G(z)\\
&{\quad}&+(z-w)\,G(w)\cdot\Lambda(w)\cdot\breve G(w)\cdot
G(z)\cdot\Lambda(z)\cdot\breve G(z)\\
&=&R^\tau(w)-R^\tau(z)
+G(w)\cdot\left(
\Lambda(z)-\Lambda(w)\right)\cdot\breve G(z)\\
&{\quad}&+(z-w)\,G(w)\cdot\Lambda(w)\cdot\breve G(w)\cdot
G(z)\cdot\Lambda(z)\cdot\breve G(z)\ ,
\end{eqnarray*}
the relation (3) is implied by (by remark 2.1, when $\tau$ has dense
range, is equivalent to) 
\begin{equation}
\Lambda(w)-\Lambda(z)=
(z-w)\,\Lambda(w)\cdot\breve G(w)\cdot
G(z)\cdot\Lambda(z)\ .
\end{equation}
Suppose now that there 
exists a (necessarily closed) operator
$$\Gamma(z):D\subseteq\X'\to\X 
$$
such that, for some open set $Z\subseteq\rho(A)$ such that $z\in Z$ iff $z^*\in Z$, one has  
$$
\forall\, z\in Z,\qquad\Gamma(z)^{-1}=\Lambda(z)\ .
$$
Then we have that (4) forces $\Gamma(z)$ to satisfy the relation
\begin{equation}
\Gamma(z)-\Gamma(w)=(z-w)\,\breve G(w)\cdot
G(z)\ ,
\end{equation}
which is equivalent to 
\begin{equation}
\forall\, \ell\in D\subseteq\X',\qquad
\der \Gamma(z)\ell=\breve G(z)\cdot
G(z)\ell\ .
\end{equation}
Regarding the identity (2), suppose that
\begin{equation}
\forall\,\ell_1,\ell_2\in D,\qquad
\ell_1(\Gamma(z^*)\ell_2)=(\,\ell_2(\Gamma(z)\ell_1)\,)^*\ .
\end{equation}
This, if $\Gamma(z)$ is densely defined, is equivalent to
$J_\X\cdot \Gamma(z^*)\subseteq \Gamma(z)'$, equality being, in the
unbounded case, stronger than (7). 
In the case where $\Gamma(z)$ has a bounded inverse given by
$\Lambda(z)$ 
as we are pretending,
(7) implies (2) which, if $\Gamma(z)$ is densely defined, is then equivalent (use e.g. 
[23, thm. 5.30, chap. III]) to
\begin{equation}
\Gamma(z)'=J_\X\cdot \Gamma(z^*)\ . 
\specialeqnum{7.1}
\end{equation}
We will therefore concentrate now on the set of maps
$$
\Gamma:\rho(A)\to \wtilde L(\X',\X)
$$
which satisfy (5) (equivalently (6)) and (7) (we are implicitly
supposing that $D$, the domain
of $\Gamma(z)$, is $z$-independent).\par
An explicit representation of the set of such maps is given by the
following
\begin{lemma} Given any $z_0\in \rho(A)$ the map 
\begin{equation}
\hat\Gamma:\rho(A)\to \wtilde
B(\X',\X)\qquad
\hat\Gamma(z):=\tau
\cdot\left(\,\frac{G(z_0)+G(z_0^*)}{2}-G(z)\,\right)
\end{equation}
satisfies (5) and (7.1).
\end{lemma}
\begin{proof} By lemma 2.1 one has 
\begin{eqnarray*}
&\quad&(z-w)\,\breve G(w)\cdot G(z)=\tau(\,G(w)-G(z)\,)\\
&=&\tau(\,G(z_0)-G(z)\,)-\tau(\,G(z_0)-G(w)\,)
\end{eqnarray*}
and so $\tau \cdot( G(z_0)-G(z))$ solves (5); by linearity 
$\hat\Gamma(z)$ is also a solution.\par
As regard (7.1) let us at first note that 
$$
J_\X\cdot\breve G(z)=\breve G(z)''\cdot J_\calH
$$
and
$$
\left(C_\calH'\cdot J_\calH y\right)(x)=\left(J_\calH y\left(C_H x\right)\right)^*
=\left(C_\calH x(y)\right)^*=\langle x,y\rangle^*=C_\calH y(x)\ .
$$
Therefore one has   
\begin{eqnarray*}
&\quad&\left(\breve G(w)\cdot G(z)\right)'=G(z)'\cdot\breve G(w)'=
\breve G(z^*)''\cdot\left(C_\calH'\right)^{-1}\cdot\breve G(w)'\\
&=&J_\X\cdot\breve G(z^*)\cdot
J_\calH^{-1}\cdot\left(C_\calH'\right)^{-1}\cdot\breve 
G(w)'=J_\X\cdot\breve G(z^*)\cdot C_\calH^{-1}\cdot\breve G(w)'\\
&=&J_\X\cdot\breve G(z^*)\cdot G(w^*)
\end{eqnarray*}
which immediately implies that $\hat\Gamma(z)$ satisfies (7.1). 
\end{proof}
\begin{remark} Lemma 2.2 shows that the set of maps
$$
\Gamma:\rho(A)\to \wtilde L(\X',\X)
$$
which satisfy (6) and (7) can be parametrized by 
$\wtilde S(\X',\X)$. Indeed, by
(6), any of these maps must differ from $\hat\Gamma(z)\in \wtilde
B(\X',\X)$ by a $z$-independent operator in $\wtilde S(\X',\X)$. 
Therefore any parametrization is of
the kind 
\begin{equation}
\Gamma_\Theta:\rho(A)\to \wtilde L(\X',\X)\qquad
\Gamma_\Theta(z)=\Theta+\Gamma(z)\,,\quad \Theta\in \wtilde S(\X',\X)\,,
\end{equation}
where $\Gamma(z)$ is some map which satisfy (6) and (7).
\end{remark}
Lemma 2.2 does not entirely solve the problem of the search of
$\Gamma(z)$ 
since $\hat\Gamma(z)$ can give rise to non-local boundary conditions 
(see remark 2.7 below); moreover $\hat
\Gamma(z)$ explicitly depends on the choice of a particular
$z_0\in\rho(A)$. However the boundedness of 
$\hat\Gamma(z)$ implies a useful
criterion for obtaining other maps $\Gamma(z)$ which satisfy (6) and
(7):
\begin{lemma} 
Suppose that $$\wtilde
\Gamma(z):D(\wtilde\Gamma) \subseteq\X'\to\X\,,\qquad z\in\rho(A)\,,
$$ 
is a family of conjugate linear, densely defined operators such that
\begin{equation}
\forall\,\ell_1,\ell_2\in D(\wtilde\Gamma)\,,\qquad
\ell_2(\wtilde\Gamma(z^*)\ell_1)=(\,\ell_1(\wtilde\Gamma(z)\ell_2)\,)^*
\end{equation}
and
\begin{eqnarray}
\forall\, \ell\in E,\, \forall\,\ell_1\in D(\wtilde\Gamma)\,,
\quad\der \ell(\wtilde\Gamma(z)\ell_1)&=&\ell
(\breve G(z)\cdot G(z)\ell_1)\\
&\equiv& \langle
G(z^*)\ell,G(z)\ell_1\rangle\ ,\nonumber
\end{eqnarray}
where $E\subseteq\X'$ is either a dense subspace or the dual of some
Schauder base in $\X$. Then $\wtilde\Gamma(z)$ is closable and its
closure satisfies (6)
and (7). 
\end{lemma}
\begin{proof} By (11) $\wtilde\Gamma(z)$ necessarily differs from (the
restriction to $D(\wtilde\Gamma)$ of) $\hat\Gamma(z)$ by a
$z$-independent, densely defined operator $\wtilde\Theta\in\wtilde
S(\X',\X)$. Being densely defined, $\wtilde\Theta$ has an adjoint and
$J_\X\cdot\wtilde\Theta\subseteq\wtilde\Theta'$. Therefore, being
$J_\X$ injective, $\wtilde\Theta$ is
closable and so, being $\hat\Gamma(z)$ bounded,
$\wtilde\Gamma(z)=\wtilde\Theta+\hat\Gamma(z)$ is closable. Denoting
by $\Theta$ the
closure of $\wtilde\Theta$, the closure of $\wtilde \Gamma(z)$ is given by 
$\Theta+\hat\Gamma(z)$, which satisfies (6) and (7) by lemma 2.2.
\end{proof}
We now state our main result:
\begin{theorem} 
Let $\Gamma_\Theta(z)$ be as in (9). Under the hypotheses 
\begin{equation}
Z_\Theta\,\neq\,\emptyset\,,
\specialeqnum{\mbox{\rm h1}}
\end{equation}
$$
Z_\Theta:=\left\{z\in\rho(A)\ :\ \exists\,
\Gamma_\Theta(z)^{-1}\in\wtilde B(\X,\X'),\ \exists\,
\Gamma_\Theta(z^*)^{-1}\in\wtilde B(\X,\X')\right\}\,,
$$
\begin{equation}
\mbox{\rm Range$\,\tau'\cap
\calH'=\left\{0\right\}$}\, ,\specialeqnum{\mbox{\rm h2}}
\end{equation}
the bounded linear operator
$$
R^\tau_\Theta(z):=R(z)+G(z)\cdot\Gamma_\Theta(z)^{-1}\cdot\breve G(z)\,,\qquad z\in
Z_\Theta\, ,
$$
is a resolvent of the self-adjoint
operator $A^\tau_\Theta$ which
coincides with $A$ on the kernel of $\tau$ and which is defined by
$$
D(A_\Theta^\tau):=\left\{\,\phi\in\calH\ :\ \phi=
\phi_z+G(z)\cdot\Gamma_\Theta(z)^{-1}\cdot\tau\,\phi_z,\quad \phi_z\in
D(A)\,\right\}\, ,
$$
$$
(-A_\Theta^\tau+z)\phi:=(-A+z)\phi_z\ .
$$
Such a definition is $z$-independent and the decomposition of $\phi$ entering
in the definition of the domain is unique. 
\end{theorem}
\begin{proof}
We have already proven that, under our hypotheses, $R_\Theta^\tau(z)$ is a
pseudo-resolvent, i.e. 
\begin{equation}
(z-w)\,R^\tau_\Theta(w)R^\tau_\Theta(z)=R^\tau_\Theta(w)-R^\tau_\Theta(z)\,.
\end{equation}
We proceed now as in the 
proof of [4,
thm. II.1.1.1]. By [23, chap. VIII, \S 1.1] $R^\tau_\Theta(z)$, being a
pseudo-resolvent, is the resolvent of a closed
operator if and only if it is
injective. Since $R^\tau_\Theta(z)\phi=0$ would imply
$$
R(z)\phi=-G(z)\cdot\Gamma_{\Theta}(z)^{-1}\cdot\breve G(z)\phi\ ,
$$
by (h2) we have $R(z)\phi=0$ (see remark 2.8) and so $\phi=0$.
\par
Since, as we have seen before, (7) implies, when $z\in Z_\Theta$,
$$
\Gamma_\Theta(z^*)^{-1}=\left(\Gamma_\Theta(z)^{-1}\right)'\cdot J_\X\,,
$$
one has
\begin{eqnarray*}
&\quad&\left(G(z)\cdot\Gamma_\Theta(z)^{-1}\cdot\breve G(z)\right)^*\\
&=&C_\calH^{-1}\cdot\left(G(z)\cdot\Gamma_\Theta(z)^{-1}\cdot\breve G(z)\right)'\cdot
C_\calH\cr
&=&C_\calH^{-1}\cdot\breve G(z)'\cdot\left(\Gamma_\Theta(z)^{-1}\right)'\cdot\breve G(z^*)''\cdot
(C_\calH')^{-1}\cdot C_\calH\cr
&=&G(z^*)\cdot\left(\Gamma_\Theta(z)'\right)^{-1}\cdot
J_\X\cdot\breve G(z^*)\cdot J_\calH^{-1}\cdot(C_\calH')^{-1}\cdot C_\calH\cr
&=&G(z^*)\cdot\Gamma_\Theta(z^*)^{-1}\cdot\breve G(z^*)\,,
\end{eqnarray*}
and so  
$$R_\Theta^\tau(z)^*=R_\Theta^\tau( z^*)\ .$$ 
This gives the denseness
of $D(A_\Theta^\tau):=\mbox{\rm Range}\,R^\tau_\Theta(z)$. Indeed
$\phi\perp 
D(A_\Theta^\tau)$, which is
equivalent to $\langle R_\Theta^\tau( z^*)\phi,\psi\rangle=0$ for all
$\psi\in\calH$, implies $\phi=0$. \par
Let us now define, on the dense domain $D(A_\Theta^\tau)$, the closed operator 
$$
A_\Theta^\tau:=R_\Theta^\tau(z)^{-1}-z
$$
which, by the resolvent identity (12), is independent of $z$;
it is self-adjoint since
$$
((A_\Theta^\tau)^*+ z^*)^{-1}=R_\Theta^\tau(z)^*=R_\Theta^\tau( 
z^*)=(A_\Theta^\tau+ z^*)^{-1}\ .
$$ 
To conclude, the uniqueness of the decomposition 
$$
\phi=
\phi_z+G(z)\cdot\Gamma_\Theta(z)^{-1}\cdot\tau\,\phi_z\,,
\qquad \phi\in D(A^\tau_\Theta)\,,
$$
is an immediate consequence of (h2).
\end{proof}
\begin{remark} 
Viewing $A$ as a bounded operator on $D(A)$ to
$\calH$, we can consider the adjoint $(-A+z^*)'$, so that
$$
(-A+z^*)'\cdot C_\calH :\calH\to D(A)'\,,\qquad (-A+z^*)'\cdot C_\calH{\,} _{\left|
D(A)\right.}= C_\calH\cdot(-A+z)
$$ 
and, by the definition of $G(z)$, 
$$
(-A+z^*)'\cdot C_\calH
\cdot G(z)=\tau'\ .
$$
Therefore, defining $Q_\phi:=\Gamma_\Theta(z)^{-1}\cdot\tau\,\phi_z$,
one has
$$
C_\calH\cdot(-A_\Theta^\tau+z)\phi=(-A+z^*)'\cdot C_\calH\phi_z
=(-A+z^*)'\cdot C_\calH\phi-\tau'Q_\phi\ ,
$$
i.e.
$$
A^\tau_\Theta\phi=C_\calH^{-1}\cdot\left(A'\cdot C_\calH\phi+\tau'Q_\phi\right)\ .
$$
Formally re-writing the last relation as
$$
A^\tau_\Theta\phi=A\phi+C_\calH^{-1}\cdot\tau'Q_\phi\ ,
$$
we can view $A^\tau_\Theta$ as a perturbation of $A$, the perturbation
being singular since, by (h2), $\tau' Q_\phi\in D(A)'\backslash \calH'$.
\end{remark} 
\begin{remark} If $\X$ is reflexive and $\Gamma_\Theta(z)$ is
densely defined, then, by (7.1), there follows 
$$
\Gamma_\Theta(z)^{-1}\in\wtilde B(\X,\X')\quad \Longrightarrow\quad
\Gamma_\Theta(z^*)^{-1}\in\wtilde B(\X,\X')\ .
$$
\end{remark}
\begin{remark} If $Z_\Theta\not=\emptyset$ then $Z_\Theta$ is
necessarily open. Indeed, by (5),
$$\Gamma_\Theta(z+h)=\Gamma_\Theta(z)+h\,\breve G(z)\cdot
G(z+h)\ ,$$ 
and so $\Gamma_\Theta(z+h)^{-1}\in\wtilde B(\X',\X)$ if $z\in Z_\Theta$
and $h$ is
sufficiently small.
\end{remark}
\begin{remark} If in the representation (9)
there exists $z_0\in\rho(A)$ such that $\Gamma(z_0)=\Gamma(z_0^*)=0$ (this is
certainly true if $\rho(A)\cap\RE\not=\emptyset$ and if one uses 
representation (8) with $z_0\in\RE$) then obviously $Z_\Theta$ is
non-empty for any invertible $\Theta\in\wtilde S(\X',\X)$. A more
significative criterion leading to (h1) will be given in proposition 2.1.
\end{remark}
\begin{remark} By the definition of $G(z)$ one has 
that (h2) is equivalent to 
\begin{equation}
\mbox{\rm Range$\,G(z)\cap
D(A)=\left\{0\right\}$}\ .\specialeqnum{\mbox{\rm h2}}
\end{equation}
\end{remark}
\begin{remark} If Kernel$\,\tau$ is dense in $\calH$ then (h2) holds true 
(it is not hard to show that the reverse implication is
false). 
Indeed the density hypothesis implies, if $Q\in\X'$, 
$$
\forall\,\psi\in\mbox{\rm
Kernel$\,\tau$},\qquad\langle\phi,\psi\rangle=Q(\tau\psi)
=\tau'Q(\psi)
\quad\Longrightarrow \quad\phi=0\ .
$$ 
This, by the definition of $G(z)$, implies
$$
R(z)\phi=
G(z)Q\quad\Longrightarrow \quad\phi=0\ ,
$$
which gives (h2). 
\end{remark}
\begin{remark} 
If in the above theorem one uses the representation
$\hat\Gamma_\Theta(z):=\Theta+\hat\Gamma(z)$ given by lemma 2.2 
one can readily check that the
domain of $A^\tau_\Theta$ is equivalently characterized in term of
``generalized boundary conditions'': $\phi\in D(A^\tau_\Theta)$ if and only if
$$
\exists\,Q_\phi\in D(\Theta)\subseteq\X'\quad\mbox{\rm such that}\quad 
\phi-\,\frac{G(z_0)+G(z_0^*)}{2}\,Q_\phi\in D(A)
$$
and
$$
\tau\left(\phi-\,\frac{G(z_0)+G(z_0^*)}{2}\,Q_\phi\right)=\Theta\, Q_\phi\ .
$$
\end{remark}
The following result states that when $\tau$ is surjective (h1) holds
true under 
relatively weak hypotheses:
\begin{proposition} Let
$\Gamma_\Theta(z)=\Theta+\Gamma(z)$ be closed, densely defined and satisfying
(5) and (7.1). If $\tau$ is surjective then 
$$
\C\backslash\RE\cup  W^-_\Theta\cup W^+_\Theta\subseteq Z_\Theta\,,
$$
where
$$
W^\pm_\Theta:=\left\{\,\lambda\in \RE\cap\rho(A)\ :\ 
\gamma(\pm\Gamma(\lambda))>-\gamma(\pm\Theta)\,\right\}\ .
$$
If $\tau$ merely has a dense range then
$$
\wtilde W^-_\Theta \cup \wtilde W^+_\Theta\subseteq Z_\Theta\,, 
$$
where
$$
\wtilde W^\pm_\Theta:=\left\{\,z\in \rho(A)\ :\ 
\frac{1}{2}\,\gamma(\pm(\Gamma(z)+\Gamma(z^*)))>-\gamma(\pm\Theta)\,\right\}\ .
$$
\end{proposition}
\begin{proof}
Writing 
$$
\Gamma(z)=\frac{1}{2}\,\left(\Gamma(z)+\Gamma(z^*)\right)+
\frac{1}{2}\,\left(\Gamma(z)-\Gamma(z^*)\right)
\equiv \Gamma_+(z)+\Gamma_-(z)\,,
$$
by (7) one has
$$
\re\, \ell(\Gamma(z)\ell)=\ell(\Gamma_+(z)\ell)\,,\qquad 
\im\, \ell(\Gamma(z)\ell)=-i\,\ell(\Gamma_-(z)\ell)\,.
$$
Thus by (5) there follows
$$
\im\, \ell(\Gamma(z)\ell)=-\frac{i}{2}\,(z-z^*)\ell(\breve G(z^*)\cdot G(z)\ell)
=\im(z)\|G(z)\ell\|_\calH^2
$$
and so, since $\Theta\in\wtilde S(\X',\X)$ implies
$\ell(\Theta\ell)\in\RE$, one has 
$$
|\ell(\Gamma_\Theta(z)\ell)|^2=
\left(\ell(\Theta \ell)+\left(\ell(\Gamma_+(z)\ell)\right)\right)^2+
\mbox{\rm Im}(z)^2\|G(z)\ell\|_\calH^4\,.
$$
Injectivity of $\Gamma_\Theta(z)$ and $\Gamma_\Theta(z)'$ 
for any $z\in \C\backslash\RE\cup
W^-_\Theta\cup W_\Theta^+$ then follows 
by injectivity of $G(z)$ (see remark 2.1), (7.1), injectivity of
$J_\X$, and the definitions of
$W^\pm_\Theta$. \par Being $\Gamma_\Theta(z)$ densely defined, one has 
$$
(\mbox{\rm Range}\,\Gamma_\Theta(z))^{\perp}=\mbox{\rm
Kernel}\,\Gamma_\Theta(z)'\,,
$$
and so injectivity of $\Gamma_\Theta(z)'$ give denseness
of the range of $\Gamma_\Theta(z)$. Being $\Gamma_\Theta(z)$ closed,
its domain is a Banach space w.r.t. the graph norm and we can apply the open mapping 
theorem to the continuous map
$\Gamma_\Theta(z):D(\Gamma_\Theta)\to\X$. Thus
to conclude the proof we need to
prove that the range of $\Gamma_\Theta(z)$ is closed. 
By [23, thm. 5.2, chap. IV]
$$
\inf\left\{\,\|G(z)\ell \|_\calH\,,\
\|\ell\|_{\X'}=1\,\right\}>0\ ,
$$
if and only if the range of $G(z)$ is closed; by the closed range
theorem (see e.g. [23, thm. 5.13, chap. IV]) the range of $G(z)$ is 
closed if and only if 
the range of $\breve G(z)$ is closed, and this is equivalent to the
range of $\tau$ being closed. Therefore, since
$$
\forall\,\ell\in D(\Gamma_\Theta)\,,\quad \|\ell\|_{\X'}=1\,,\qquad
\|\Gamma_\Theta(z)\ell \|_\X\ge |\ell(\Gamma_\Theta(z)\ell)|\,,
$$
when either $z\in\C\backslash\RE\cup W^-_\Theta\cup W^+_\Theta$ if $\tau$
is surjective, or when $z\in\wtilde W^-_\Theta\cup \wtilde W^+_\Theta$ if $\tau$
has a dense range, one has
$$
\inf\left\{\,\|\Gamma_\Theta(z)\ell \|_\X\,,\
\|\ell\|_{\X'}=1\,\right\}>0\ ,
$$
and so, since $\Gamma_\Theta(z)$ is closed, it has a closed range by 
[23, thm. 5.2, chap. IV].
\end{proof}
Since $Z_\Theta\subseteq\rho(A^\tau_\Theta)$, the above proposition
immediately implies a semi-boundedness criterion for the extensions
$A_\Theta^\tau$:
\begin{corollary} Let $-A$ be bounded from below and suppose
that there exist $\lambda_0 \in\rho(A)\cap\RE$ and $\theta_0\in\RE$ such that 
$$
\forall\,\lambda\ge\lambda_0\qquad
\gamma(\Gamma(\lambda))>-\theta_0\ .
$$
Then
$$\inf\sigma(-A_\Theta^\tau)\ge-\lambda_0
$$  
for any $\Theta\in\wtilde S(\X',\X)$ such that $\gamma(\Theta)\ge
\theta_0$.
\end{corollary}
\begin{remark} 
By the proposition above, if $\X=\,$Range$\,\tau$ is
finite-dimensional and $\Gamma_\Theta(z)$ is
everywhere defined, then (h1) is satisfied with at least 
$\C\backslash\RE\subseteq Z_\Theta$.
\end{remark}
\begin{remark} By the proposition above, since $\hat \Gamma(z)$ is
bounded, 
if one uses the 
representation $\hat \Gamma_\Theta(z)$, with
$\Theta\in \wtilde S(\X',\X)$ closed, densely defined and such that $J_\X\cdot\Theta=\Theta'$, then (h1) is
satisfied (with at least $\C\backslash\RE\subseteq Z_\Theta$) when 
$\tau$ is surjective.
\end{remark}
\begin{remark} 
If $\X$ is a Hilbert space 
(with scalar product $\langle\cdot,\cdot\rangle$) we can of course use the map 
$C_\X$ to identify $\X$ with $\X'$ and 
re-define $G(z)$ as
$$
G(z):=C_\calH^{-1}\cdot\breve G(z^*)'\cdot C_\X:\X\to\calH\ .
$$
The statements in the above theorem remain unchanged taking
$$
\Gamma_\Theta:\rho(A)\to L(\X,\X)\,,\qquad
\Gamma_\Theta(z)=\Theta+\Gamma(z)
$$
with $\Theta$ such that  
$$
\forall\, x,y\in D(\Theta),\qquad \langle \Theta x,y\rangle=
\langle x,\Theta y\rangle
$$
and $\Gamma(z)$ satisfying (6) and
$$
\forall\, x,y\in D(\Gamma),\qquad \langle \Gamma(z)x,y\rangle=\langle
x,\Gamma(z^*)y\rangle\ .
$$
\end{remark}
\begin{remark} When $\X$ is a Hilbert space, by theorem 2.1, since
$G(z)$ and $\breve G(z)$ are bounded, we have that
$R^\tau_\Theta(z)-R(z)$ is a trace class operator on $\calH$ if and only if
$(\Theta+\Gamma(z))^{-1}$ is a trace class operator on $\X$ (see
e.g. [23, \S 1.3, chap. X]). This
information can be used (proceeding along the same lines as in \cite{[BT]})
to infer from $\sigma(A)$ some properties of $\sigma(A^\tau_\Theta)$. 
\end{remark}
When $\X$ is a Hilbert space one can give, besides the one appearing
in lemma 2.2, another criterion for obtaining the map
$\Gamma_\Theta$. Indeed one has the following
\begin{lemma} 
Suppose that there exists a densely
defined sesquilinear form 
$\wtilde\E(z)$, $z\in\rho(A)$, with $z$-independent domain
$D(\wtilde\E)\times 
D(\wtilde\E)$, such
that 
\begin{equation}
\forall\,x,y\in
D(\wtilde\E),\qquad\wtilde\E(z^*)(x,y)=(\,\wtilde\E(z)(y,x)\,)^*\ ,
\end{equation}
\begin{equation}
\forall\,x,y\in D(\wtilde\E),\qquad\der\wtilde\E(z)(x,y)=\langle
G(z^*)x,G(z)y\rangle\ ,
\end{equation}
and such that there exist $z_0\in\rho(A)$, $M\in\RE$ for which $\wtilde\E(z_0)$
is closable and   
\begin{equation}
\forall\,x\in
D(\wtilde\E),\qquad\mbox{\rm Re}\left(\wtilde\E(z_0)(x,x)\right)
\ge M\,\langle x,x\rangle\ .
\end{equation}
Then $\wtilde\E(z)$
is closable for any $z\in\rho(A)$ and, denoting its closure by $\E(z)$, 
there exists a densely 
defined, closed linear operator
$\Gamma(z)$ with $z$-independent domain $D(\Gamma)$, defined by
$$
\forall\,x\in D(\E),\,\forall\,y\in D(\Gamma),\qquad\E(z)(x,y)
=\langle x,\Gamma(z)y\rangle\ ,$$ 
satisfying (6) and the Hilbert space analogue of (7.1), i.e.
$$\Gamma(z)^*
=\Gamma(z^*)\ .$$
\end{lemma}  
\begin{proof} By our hypotheses $\wtilde\E(z)$ necessarily differs from (the
restriction to $D(\wtilde\E)\times D(\wtilde\E)$ of) the bounded
sesquilinear form associated to $\hat\Gamma(z)$ by a
$z$-independent Hermitean form $\wtilde\mathcal Q$. Therefore 
$$\wtilde\Q(x,y)=\wtilde\E(z_0)(x,y)-\langle x,\hat\Gamma(z_0)y\rangle$$
is a semi-bounded, densely defined, closable Hermitean form. If $\Theta$ denotes the
unique semi-bounded self-adjoint operator corresponding to the closure of
$\wtilde\Q$ (see [23, thm. 2.6,
chap. VI] for the existence of $\Theta$), then 
the operator $\Gamma(z):=\Theta+\hat\Gamma(z)$ gives the thesis. 
\end{proof}
\begin{remark}
If $\wtilde\Gamma(z)$ in lemma 2.3, besides satisfying
(10) and (11), is bounded from below
in the sense 
of (15), i.e. if there exist $z_0\in\rho(A)$, $M\in\RE$ such that
$$
\forall\,x\in D(\wtilde\Gamma(z_0))\,,
\qquad\re(\langle x,\wtilde\Gamma(z_0)x\rangle)\ge 
M\,\langle x,x\rangle\ ,
$$
then, by using both lemma 2.3 and lemma 2.4, it is closable and
its closure satisfies (5) and (7.1). This is
nothing but a variation of Friedrichs extension theorem.  
\end{remark}
\begin{remark}
The operator $\Gamma(z)$ given by lemma 2.4 satisfies
$$
\frac{1}{2}\,\gamma(\Gamma(z_0)+\Gamma(z_0^*))\ge M
$$
and so, when $\tau$ has dense range, by proposition 2.1 one has that
(h1) holds true 
for $\Gamma_\Theta(z)$, where $\Theta$ is any $\Gamma(z)$-bounded 
(see [23, thm. 1.1, chap. IV) self-adjoint operator such that
$\gamma(\Theta)>-M$. If the 
constant $M$ can be made arbitrarily large by
letting $|z_0|\uparrow\infty$, then (h1) is satisfied with any bounded
from below self-adjoint operator $\Theta$.
\end{remark}
\section{APPLICATIONS}
\begin{example} {\it The $\calH_{-2}$ -construction.} Let $\X=\C$, $\varphi\in D(A)'\backslash\left\{0\right\}$ and put $\tau=\varphi$. Defining 
$$
\wtilde
R(z):=C_\calH^{-1}\cdot R(z^*)'
\in B(D(A)',\calH)$$
one has then
$$
\breve G(z):\calH\to\C\,,\qquad\breve G(z)\phi=\langle\wtilde R(z^*)\varphi,\phi\rangle
$$
and
$$
G(z):\C\to\calH\,,\qquad G(z)\zeta=\zeta\wtilde R(z)\varphi\ .
$$
The hypothesis (h2) is equivalent to the request
$$
\varphi\notin \calH'\ ,
$$
whereas hypothesis (h1) is always satisfied with at least $\C\backslash\RE\subseteq
Z_\Theta$ since $\X$ is finite dimensional
(see remark 2.11). Then the self-adjoint operator 
$A^\varphi_\alpha$ has resolvent
$$
(-A^\varphi_\alpha+z)=(-A+z)^{-1}+\Gamma_\alpha(z)^{-1}\wtilde
R(z)\varphi\otimes
\wtilde R(z^*)\varphi\ ,
$$
where (by lemma 2.2)
$$
\Gamma_\alpha(z)
=\alpha+\varphi\left(\frac{\wtilde R(z_0)\varphi+\wtilde
R(z_0^*)\varphi}{2}-
\wtilde R(z)\varphi\right)\,,\qquad\alpha\in\RE\ .
$$
This coincides with the $``\,\calH_{-2}$ -construction'' given in \cite{[KS]}
(there only the case $-A\ge 0$, $z_0=1$ was considered). For a similar
construction also see \cite{[AK]} and references therein.
\end{example}
\begin{example} {\it A variation on the Birman-Kre\u\i n-Vishik
theory.} Let $A$ be a strictly positive self-adjoint operator, so that $0\in
\rho(-A)$, and let $\tau: D(A)\to\X$ satisfy (h2). By remark 2.7  and
theorem 2.1, for any $\Theta\in\wtilde S(\X',\X)$ which has a bounded inverse, we can
define  the (strictly positive when $\Theta$ is positive,
i.e. $\gamma(\Theta)\ge 0$) 
self-adjoint opertor $A^\tau_\Theta$ by
$$
({A^\tau_\Theta})^{-1}=A^{-1}+G\cdot\Theta^{-1}\cdot\breve G\ ,
$$
where $G:=G(0)$ and $\breve G:=\breve G(0)$. Moreover one has 
$$
D(A_\Theta^\tau):=\left\{\,\phi\in\calH\ :\ \phi=
\phi_0+G Q_\phi,\ \phi_0\in
D(A),\ \tau\phi_0=\Theta\,Q_\phi\,\right\}\,,
$$
$$
A_\Theta^\tau\phi=A\phi_0\ .
$$
This gives a variation of the 
Birman-Kre\u\i n-Vishik approach which comprises the result given in
\cite{[KKO]}. In particular [22, example 4.1] can be obtained
by taking $\calH=L^2(\Omega)$, $A=-\Delta_\Omega+\lambda$, $\lambda>0$,
$\Omega=(0,\pi)\times\RE^2$, $D(\Delta_\Omega)=H^2_0(\Omega)$,
$\tau:H^2_0(\Omega)\to L^2(0,\pi)$
the evaluation along the segment $\left\{(x,0,0), x\in (0,\pi)\right\}$,
$\Theta=
-\Delta_{(0,\pi)}$, $D(\Theta)=H^2_0(0,\pi)$; [22, example 4.2]
corresponds to $\calH=L^2(\RE^3)$, $A=-\Delta+\lambda$, $\lambda>0$, 
$D(\Delta)=H^2(\RE^3)$, whereas $\tau$ and $\Theta$ are the same as before.
\end{example}
\begin{example} {\it Finitely many point interactions in three
dimensions.} We take $\calH=L^2(\RE^3)$, $A=\Delta$, $D(A)=H^2(\RE^3)\subset
C_b(\RE^3)$. Considering then a finite set $Y\subset\RE^3$, $\#Y=n$, 
we take as the linear operator $\tau$ the linear continuous surjective map
$$
\tau_Y:H^2(\RE^3)\to\C^n\qquad\tau_Y\phi:=\{\phi(y)\}_{y\in Y}\ .
$$
Then one has 
$$
\breve G(z):\LD\to\C^n\,,\qquad\breve G(z)\phi=\{\G_z*\phi\,(y)\}_{y\in Y}\ ,
$$
where $$\G_z(x)=\frac{e^{-\sqrt z|x|}}{4\pi|x|}\,,\qquad \re\sqrt
z>0\,,\qquad\G_z^y(x):=\G_z(x-y)\ ,$$ 
and
$$
G(z):\C^n\to\LD\,,\qquad G(z)\zeta=\sum_{y\in Y}\zeta_y\G_z^y
\equiv \G_z*\sum_{y\in Y}\zeta_y\delta_y\ .
$$
A straightforward calculation then gives 
\begin{eqnarray*}
&\quad&\left(\breve G(z)\cdot G(z)\zeta\right)_y
=\sum_{\tilde y\in Y}\zeta_{\tilde y}
\langle\G^{\tilde y}_{z^*},\G^y_z\rangle\cr
&=&\zeta_y\,\frac{1}{(2\pi)^3}
\int_{\RE^3} dk\,
\frac{1}{(|k|^2+z)^2}
+\sum_{\tilde y\not=y}\zeta_{\tilde y}\,\frac{1}{(2\pi)^3}
\int_{\RE^3} dk\,
\frac{e^{-ik\cdot(\tilde y-y)}}{(|k|^2+z)^2}\cr
&=&\zeta_y\,\frac{1}{2\pi^2}
\int_0^\infty dr\,
\frac{r^2}{(r^2+z)^2}
+\sum_{\tilde y\not=y}\zeta_{\tilde y}\,\frac{1}{2\pi^2|\tilde y-y|}
\int_0^\infty dr\,
\frac{r\sin(r|\tilde y-y|)}{(r^2+z)^2}\cr
&=&\zeta_y\,\frac{1}{8\pi\sqrt z}+\sum_{\tilde
y\not=y}\zeta_{\tilde y}\,\frac{e^{-\sqrt z\,|\tilde y-y|}}{8\pi\sqrt z}
=\der\left(\zeta_y\,\frac{\sqrt z}{4\pi}-\sum_{\tilde
y\not=y}\zeta_{\tilde y}\,\G_z^{\tilde yy}\,\right)
\ ,
\end{eqnarray*}
where $\G_z^{\tilde yy}:=\G_z(\tilde y-y)$, $\tilde y\not=y$. Defining 
$$
\wtilde\G_z:\C^n\to\C^n\qquad(\wtilde \G_z\zeta)_y:=\sum_{\tilde
y\not=y}\zeta_{\tilde y}\,\G_z^{\tilde yy}\ ,
$$
one can take as $\Gamma_\Theta(z)$ the linear operator
$$
\Gamma_\Theta(z)=\Theta+\frac{\sqrt z}{4\pi}-\wtilde\G_z\ ,
$$
where $\Theta$ is any Hermitean $n\times n$ matrix.\par 
Hypothesis (h1) is satisfied with at least $\C\backslash\RE\subseteq
Z_\Theta$ since $\X$ is finite dimensional
(see remark 2.11) and hypotheses (h2) is satisfied since $\G_z^y\notin
H^2(\RE^3)$ for any $y\in Y$. In conclusion on can define the
self-adjoint operator 
$\Delta_{\Theta}^Y$ with resolvent given by
$$
(-\Delta_{\Theta}^Y+z)^{-1}=(-\Delta+z)^{-1}+
\sum_{y,\tilde y\in Y}\left(\Theta+\frac{\sqrt z}{4\pi}-\wtilde\G_z
\right)^{-1}_{y\tilde y}\G_z^y\otimes\G_{z^*}^{\tilde y}\ .
$$
This coincides with the operator
constructed in [4, \S II.1.1].
\end{example}
\begin{example} {\it Infinitely many point
interactions in 
three
dimensions.} We take $\calH=L^2(\RE^3)$, $A=\Delta$, $D(A)=H^2(\RE^3)\subset
C_b(\RE^3)$. Considering then an infinite and countable set
$Y\subset\RE^3$ such that
\begin{equation}
\inf_{y\not=\tilde y\ y,\tilde y\in Y}\, |y-\tilde y|=d>0\ ,
\end{equation}
we take as the linear operator $\tau$ the linear map
$\tau_Y\phi:=\{\phi(y)\}_{y\in Y}$. The hypothesis (16) ensures its
surjectivity and (see
[4, page 172])
$$
\tau_Y:H^2(\RE^3)\to\ell_2(Y)\,,\qquad\tau_Y\in B(H^2(\RE^3),\ell_2(Y))
$$ 
Proceeding as in the previous example one has then
$$
\breve G(z):\LD\to\ell_2(Y)\,,\qquad\breve G(z)\phi=\{\G_z*\phi\,(y)\}_{y\in Y}\ ,
$$
and
$$
G(z):\ell_2(Y)\to\LD
$$
is the unique bounded linear operator which, on the dense subspace 
$$\ell_0(Y):=\left\{\zeta\in\ell_2(Y)\ :\ \#\,\mbox{\rm
supp($\zeta$)}<+\infty\right\}\, ,$$ is defined by 
$$
G(z)\zeta=\sum_{y\in Y}\zeta_y\G_z^y\ ,
$$ 
i.e.
$$\forall\,\zeta\in\ell_2(Y),\qquad G(z)\zeta=\G_z*\tau'_Y(\zeta)\ ,$$
where $\tau'_Y(\zeta)\in H^{-2}(\RE^3)$ is the signed Radon measure
defined by $$\tau'_Y(\zeta)(\phi)=\langle\zeta^*,\tau_Y\phi\rangle\, .$$ 
Taking $\zeta\in\ell_0(Y)$ one then obtains, proceeding as in example 3.2, 
\begin{eqnarray*}
\left(\breve G(z)\cdot G(z)\zeta\right)_y
&=&\zeta_y\,\frac{1}{8\pi\sqrt z}+\sum_{\tilde
y\not=y}\zeta_{\tilde y}\,\frac{e^{-\sqrt z\,|\tilde y-y|}}{8\pi\sqrt z}\\
&=&\der\left(\zeta_y\,\frac{\sqrt z}{4\pi}-\sum_{\tilde
y\not=y}\zeta_{\tilde y}\,\G_z^{\tilde yy}\,\right)
\ .
\end{eqnarray*}
Posing
$$
\wtilde\G_z:\ell_0(Y)\to\ell_2(Y)\,,\qquad
(\wtilde\G_z\zeta)_y:=\sum_{\tilde
y\not=y}\zeta_{\tilde y}\,\G_z^{\tilde yy}\ ,
$$
the operator
$$
\wtilde \Gamma(z):=\frac{\sqrt z}{4\pi}-\wtilde\G_z
$$
satisfies (10) and (11). Therefore, by lemma 2.3
(with $E$ the canonical basis of $\ell_2(Y)$ and
$D(\wtilde\Gamma)=\ell_0(Y)$), $\wtilde\G_z$ is closable and, denoting
its closure by the same symbol, the closed and densely
defined operator
$$
\Gamma(z):=\frac{\sqrt z}{4\pi}-\wtilde\G_z\ ,
$$
satisfies (5) and (7). Since $\Gamma(z)+\Gamma(z^*)$ is bounded from below if
$\im(z)$ is sufficiently large (see [4, page 171]), by lemma 2.4 it satisfies (7.1).
Therefore, considering then $\Theta+\Gamma(z)$, where 
$\Theta$ is any $\Gamma(z)$-bounded (see [23,
thm. 1.1. chap. IV]) self-adjoint operator on $\ell_2(Y)$, 
(h1) is satisfied by proposition 2.1, whereas (h2) is equivalent to 
$\tau'_Y(\zeta)\notin L^2(\RE^3)$ for any
$\zeta\not=0$, which is always true since the support of $\tau'_Y(\zeta)$
is the null set $Y$. So, by theorem 2.1, one can define the self-adjoint
operator 
$\Delta_{\Theta}^Y$ with resolvent given by
$$
(-\Delta_{\Theta}^Y+z)^{-1}=(-\Delta+z)^{-1}+
\sum_{y,\tilde y\in Y}\left(\Theta+\frac{\sqrt z}{4\pi}-\wtilde\G_z
\right)^{-1}_{y\tilde y}\G_z^y\otimes\G_{z^*}^{\tilde y}
\ .
$$
This coincides with the operator constructed 
(by an approximation
method) in [4, \S III.1.1]. 
\end{example}
\begin{example} {\it Singular perturbations of the Laplacian 
supported by regular curves.} We take $\calH=L^2(\RE^n)$, $A=\Delta$, $D(A)=H^2(\RE^n)$, $n=3$ or $n=4$. 
Consider then a $C^2$ curve $\gamma:I\subseteq\RE\to\RE^n$ such that
$C:=\gamma(I)$ is a one-dimensional embedded submanifold 
$C\subset\RE^n$ which, when unbounded, is, outside
some compact set, 
globally diffeomorphic to a straight line (these hypotheses on
$\gamma$ will be
weakened in the next example). We will suppose $C$ to be parametrized
in such a way that $|\dot\gamma|=1$.\par
We take as linear
operator $\tau$ the unique linear map 
$$
\tau_\gamma:H^2(\RE^n)\to L^2(I)\,,\qquad\tau_\gamma\in B(H^2(\RE^n),L^2(I))
$$
such that
$$
\forall\phi\in C^\infty_0(\RE^n)\,,\qquad\tau_\gamma\phi(s):=\phi(\gamma(s))\ .
$$
The existence of such a map is given by combining the results in [8,
\S10] (straight line) with the ones in [8, \S 24] (compact manifold).
By [8, \S 25] we have that 
$$\mbox{\rm Range}\,\tau_\gamma=H^s(I),\quad
\tau_\gamma\in B(H^2(\RE^d),H^s(I)),\qquad  
s=2-\,\frac{n-1}{2}
$$
and so we could take $\X=H^s(I)$. However, in order to make clearer
the connections with the existing literature, we prefer 
to work with $\X=L^2(I)$ even if with this choice
$\tau_\gamma$ is not surjective (but has a dense range). 
\smallskip\p
\paragraph{The case n=3}
One has, proceeding similarly to examples 3 and 4, 
$$
\breve G(z):\ld3\to L^2(I)\,,\qquad\breve G(z)\phi
=\tau_\gamma(\G_z*\phi)
$$
and
$$
G(z):L^2(I)\to\ld3\,,\qquad G(z)f=\G_{z}*\tau'_\gamma(f)\ ,
$$ 
where $\tau'_\gamma(f)\in H^{-2}(\RE^3)$ is the signed Radon measure
defined by $$\tau'_\gamma(f)(\phi)=\langle f^*,\tau_\gamma\phi\rangle\,.$$ 
By Fourier transform one has equivalently
$$
\F\cdot G(z)f(k)=\frac{1}{(2\pi)^{{3}/{2}}}\ \frac{1}{|k|^2+z}\,
\int_I ds\,f(s)\,e^{-ik\cdot\gamma(s)}\,,\qquad f\in L^1(I)\cap L^2(I)
\ ,
$$
so that, for any $f_1,f_2\in L^1(I)\cap L^2(I)$ one obtains 
\begin{eqnarray}
&\quad&(z-w)\langle f_1,\breve G(w)\cdot G(z)f_2\rangle\nonumber\\
&=&(z-w)\,\langle G(w^*)f_1, G(z)f_2\rangle\nonumber\\
&=&\frac{(z-w)}{(2\pi)^3}\,\int_{I^2} dt\,ds\,{f^*_1(t)\,f_2(s)}
\int_{\RE^3} dk\,
\frac{e^{-ik\cdot(\gamma(t)-\gamma(s))}}{(|k|^2+w)\,(|k|^2+z)}\nonumber\\
&=&\frac{(z-w)}{2\pi^2}\,\int_{I^2} dt\,ds
\,\frac{f_1^*(t)\,f_2(s)}{|\gamma(t)-\gamma(s)|}
\int_0^\infty dr\,
\frac{r\sin(r|\gamma(t)-\gamma(s)|)}{(r^2+w)\,(r^2+z)}\nonumber\\
&=&\int_{I^2} dt\,ds\,f_1^*(t)\,f_2(s)
\ \frac{e^{-\sqrt w\,|\gamma(t)-\gamma(s)|}-
e^{-\sqrt z\,|\gamma(t)-\gamma(s)|}}{4\pi\,|\gamma(t)-\gamma(s)|}
\ .
\end{eqnarray}
Suppose now that, in the case $I$ is not compact, 
\begin{equation}
\exists\,\lambda_0>0\ :\ \forall\,\lambda\ge \lambda_0,\qquad
\sup_{t\in I}\int_{I}ds\, e^{-\lambda|\gamma(t)-\gamma(s)|}<+\infty\ .
\end{equation}
By (17) one can then define a linear operator ${\breve\Gamma_\epsilon}(z):L_0^2(I)\to
L^2(I)$, $\epsilon>0$, satisfying (5) and (7), by
$$
\breve\Gamma_\epsilon(z)f(t):=
\int_{I} ds\,f(s)\left(
\,\frac{\chi_\epsilon(t,s)}{4\pi\,|t-s|}-\frac{e^{-\sqrt
z\,|\gamma(t)-\gamma(s)|}}{4\pi\,|\gamma(t)-\gamma(s)|}\,\right)\ ,
$$
where $\chi_\epsilon(t,s):=\chi_{[0,\epsilon]}(|t-s|)$ and $$L^2_0(I):=\left\{f\in
L^2(I)\, :\, \mbox{\rm $f$ has compact support}\right\}\,.$$
When $f\in C_0^1(I)$ one can then re-write
$\breve\Gamma_\epsilon(z)f$ as
\begin{eqnarray*}
\breve\Gamma_\epsilon(z)f(t)&=&\int_{I} ds\,(\,f(t)-f(s)\,)
\,\G_z(\gamma(t)-\gamma(s))\cr
&\quad&+f(t)\int_{I} ds
\,\frac{\chi_\epsilon(t,s)}{4\pi\,|t-s|}-\frac{e^{-\sqrt
z\,|\gamma(t)-\gamma(s)|}}{4\pi\,|\gamma(t)-\gamma(s)|}\cr
&\quad&-\int_{I} ds\,\chi_\epsilon(t,s)
\,\frac{f(t)-f(s)}{4\pi\,|t-s|}
\ .
\end{eqnarray*}
The second term has,
as a function of the parameter $\epsilon>0$, a derivative given by 
$(2\pi\epsilon)^{-1}f(t)$, and the
last term is $z$-independent. Therefore the operator 
$\wtilde\Gamma(z):C_0^1(I)\to L^2(I)$,
\begin{eqnarray}
\wtilde\Gamma(z)f(t):=
\int_{I} ds\,(\,f(t)-f(s)\,)
\,\G_z(\gamma(t)-\gamma(s))&\nonumber\\
+f(t)\left(\,\frac{1}{2\pi}\,\log\left(\epsilon^{-1}\right)+\int_{I} ds
\,\frac{\chi_\epsilon(t,s)}{4\pi\,|t-s|}-\frac{e^{-\sqrt
z\,|\gamma(t)-\gamma(s)|}}{4\pi\,|\gamma(t)-\gamma(s)|}
\,\right)&
\end{eqnarray}
is $\epsilon$-independent and satisfies (10) and (11) with $E=L^2(I)$ and
$D(\wtilde\Gamma)=C_0^1(I)$. Moreover, by lemma 2.3, it is closable and its
closure $\Gamma(z)$ satisfies (6) and (7). Since
$\Gamma(z)+\Gamma(z^*)$ is bounded 
from below if
$\im(z)$ is sufficiently large (this is a conseguence of (18)), 
by remark 2.15 it satisfies (7.1). Moreover (see [38, lemma 1]) such a bound can be made
arbitrarily large by letting $|z|\uparrow\infty$.
Therefore, considering then $\Theta+\Gamma(z)$, where 
$\Theta$ is any $\Gamma(z)$-bounded self-adjoint operator on $L^2(I)$, 
by remark 2.16 and
proposition 2.1, (h1) is satisfied when $\Theta$ is bounded from below,
whereas (h2) is satisfied since $\tau'_\gamma (f)\notin L^2(\RE^3)$ for
any $f\not=0$, being the support of $\tau'_\gamma (f)$ given by the null
set $C$. \par  
The corresponding self-adjoint family given by theorem 2.1 has resolvents
$$
(-\Delta_\Theta^\gamma+z)^{-1}\phi=(-\Delta+z)^{-1}\phi
+\G_z*\tau'_\gamma\left(\,(\,\Theta+
\Gamma(z)\,)^{-1}\cdot\tau_\gamma(\G_z*\phi)\,\right)\ .
$$
These give singular perturbations of the Laplacian of the same kind
obtained (by a quadratic form approach) in \cite{[T]}. 
\smallskip\p
\paragraph{The case n=4}
Proceeding as in the case $n=3$ one obtains
$$
\breve G(z):\ld4\to L^2(I)\,,\qquad\breve G(z)\phi
=\tau_\gamma(\K_z*\phi)\ ,
$$
$$
G(z)f=\K_{z}*\tau'_\gamma(f),\quad
\F\K_z(k):=\frac{1}{|k|^2+z},\quad k\in\RE^4
$$
and, for any $f_1,f_2\in L^1(I)\cap L^2(I)$,
\begin{eqnarray*}
&\quad&(z-w)\langle f_1,\breve G(w)\cdot G(z)f_2\rangle
=(z-w)\,\langle G(w^*)f_1, G(z)f_2\rangle\cr
&=&\int_{I^2} dt\,ds\,f_1^*(t)\,f_2(s)
\,\left(\,\K_w(\gamma(t)-\gamma(s))-\K_z(\gamma(t)-\gamma(s))\,\right)
\ .
\end{eqnarray*}
Since 
\begin{eqnarray*}
|\K_z(x)|&=&\frac{1}{4\pi^2|x|^2}\,(\,1+o\left(|x|\right)\,)\,,\qquad |x|\ll
1\ ,\cr
|\K_z(x)|&=&\frac{1}{2\,(2\pi)^{3/2}|x|^{3/2}}\,e^{-2{{\mathrm {Re}}}\sqrt
z\,|x|}\,
(\,1+o\left(1/|x|\right)\,)\,,\qquad |x|\gg
1\ ,
\end{eqnarray*}
when $\gamma$ satisfies (18) the linear operator 
$\breve\Gamma_\epsilon(z):L_0^2(I)\to
L^2(I)$ 
$$
\breve\Gamma_\epsilon(z)f(t)=
\int_{I} ds\,f(s)
\,\left(\,\frac{\chi_\epsilon(t,s)}{4\pi^2\,|t-s|^2}-\K_z(\gamma(t)-\gamma(s))\,\right)
$$
is well defined and 
satisfies (10) and (11) with $E=L^2(I)$ and $D(\breve\Gamma)=L_0^2(I)$. \par
In four dimensions, due to the stronger (w.r.t. $\G_z$) singularity at
the origin of $\K_z$, it is no longer possible to perform the
calculations leading to the analogue of the operator $\wtilde
\Gamma(z)$, and one is forced to use sesquilinear forms and to try then
to apply lemma 2.4. Defining for brevity
$$
k_{(\epsilon)}(t,s):=\frac{\chi_\epsilon(t,s)}{4\pi^2\,|t-s|^2}\,,\qquad
k_z(t,s):=\K_z(\gamma(t)-\gamma(s))\ ,
$$ one can re-write $\langle f_1,\breve\Gamma_\epsilon(z)f_2\rangle$, when 
$f_1,f_2\in C_0^1(I)$, as
\begin{eqnarray*}
&\quad&\langle f_1,\breve\Gamma_\epsilon(z)f_2\rangle=\int_{I^2}dt\,
ds\,f_1^*(t) f_2(s)\,(\,k_{(\epsilon)}(t,s)-k_z(t,s)\,)\cr
&=&\int_{I^2}dt\,
ds\,(\,f_1^*(t) f_2(s)-f_1^*(t) f_2(t)+f_1^*(s) f_2(s)\,)\,(\,k_{(\epsilon)}(t,s)-k_z(t,s)\,)\cr
&=&\frac{1}{2}\,\int_{I^2}dt\,
ds\,(\,f_1^*(t)-
f^*_1(s)\,)\,(\,f_2(t)-f_2(s)\,)\,(\,k_z(t,s)-k_{(\epsilon)}(t,s)\,)\cr
&\quad&+\int_{I^2}dt\,
ds\,f_1^*(t)f_2(t)\,(\,k_{(\epsilon)}(t,s)-k_z(t,s)\,)\cr
&=&\frac{1}{2}\,\int_{I^2}dt\,
ds\,(\,f_1^*(t)-
f^*_1(s)\,)\,(\,f_2(t)-f_2(s)\,)\,k_z(t,s)\cr
&\quad&+\int_{I^2}dt\,
ds\,f_1^*(t)f_2(t)\,(k_{(\epsilon)}(t,s)-k_z(t,s))\cr
&\quad&-\frac{1}{2}\,\int_{I^2}dt\,
ds\,(\,f_1^*(t)-
f^*_1(s)\,)\,(\,f_2(t)-f_2(s)\,)\,k_{(\epsilon)}(t,s)
\ .
\end{eqnarray*}
Similarly to the three dimensional case the second term has,
as a function of the parameter $\epsilon>0$, a derivative given by 
$(2\pi^2\epsilon^2)^{-1}\smallint_I dt\, f^*_1(t)f_2(t)$, and the
last term is $z$-independent. Therefore the
sesquilinear form
$$\wtilde\E(z):C_0^1(I)\times C_0^1(I) \to\C\ ,$$
\begin{eqnarray}
\wtilde\E(z)(f_1,f_2):=
\frac{1}{2}\,\int_{I^2}dt\,
ds\,(\,f_1^*(t)-
f^*_1(s)\,)\,(\,f_2(t)-f_2(s)\,)\,k_z(t,s)&\quad&\nonumber\\
+\int_{I}dt\,f_1^*(t)f_2(t)\,\left(\,\frac{1}{2\pi^2\epsilon}+
\int_I ds\,(k_{(\epsilon)}(t,s)-k_z(t,s))\,\right)&\quad&
\end{eqnarray}
is $\epsilon$-independent and satisfies (13) and (14). 
It is straightforward to check its closability (see the proof of proposition
2 in \cite{[T]} if you get stuck), whereas (15) is a
consequence of (18). Moreover, proceeding as in the case $n=3$, 
the bound in (15) can be made
arbitrarily large by letting $|z|\uparrow\infty$. As (h2) is verified by
the same 
argument as in
the case $n=3$, by lemma 2.4, remark 2.16, proposition 2.1 and 
theorem 2.1, one has a self-adjoint family of
self-adjoint operators with resolvents
$$
(-\Delta_\Theta^\gamma+z)^{-1}\phi=(-\Delta+z)^{-1}\phi
+\K_z*\tau'_\gamma\left(\,(\Theta+\Gamma(z)\,)^{-1}\cdot\tau_\gamma(\K_z*\phi)\,\right)\ ,
$$
where $\Gamma(z)$ is the operator corresponding to the closure of
$\wtilde\E(z)$ and $\Theta$ is any bounded from below self-adjoint
operator on $L^2(I)$. 
This gives singular perturbations of the Laplacian of the same kind obtained in
\cite{[Sh]}. \par
\end{example}
\begin{example} {\it Singular perturbations given by $d$-sets and $d$-measures.} 
A Borel set $F\subset\RE^n$ is called a $d$-set, $d\in(0,n]$, if (see [21,
Chap. II]) there exists a Borel
measure $\mu$ in $\RE^n$ such that supp$\,(\mu)=F$ and
\begin{equation}
\exists\, c_1,\,c_2>0\,:\, \forall x\in F,\ \forall r\in(0,1),\ 
c_1r^d\le\mu(B_r(x)\cap F)\le c_2r^d,
\end{equation}
where $B_r(x)$ is the ball of radius $r$ centered at the point
$x$. By [21, chap. II, thm. 1], once $F$ is a $d$-set, $\mu_F$, the
$d$-dimensional Haurdorff measure restricted to $F$, always satisfies
(21) and so $F$ has
Hausdorff dimension $d$ in the neighbourhood of any of its
points. From the definition there also follows that a finite union of
$d$-sets which intersect on a set of zero $d$-dimensional Hausdorff
measure is a $d$-set.
Examples of $d$-sets are $d$-dimensional Lipschitz manifolds (use examples 2.1 and 2.4
in \cite{[JW1]}) and (when $d$ is not an integer) self-similar fractals of Hausdorff
dimension $d$ (see [21, chap. II, example 2], [39, thm. 4.7]).   \par
Denoting by $j_F:F\to\RE^n$ the restriction to the $d$-set $F$ of the
identity map, 
we take as the linear operator $\tau$ the unique continuous map ($0<n-d<2s$) 
$$
\tau_F:H^s(\RE^n)\to L^2(F)\,,\qquad\tau_F\in B(H^s(\RE^n),L^2(F))
$$
such that
$$
\forall\phi\in C^\infty_0(\RE^n)\,,\qquad\tau_F\phi(x):=\phi(j_F(x))\ .
$$
Here $L^2(F)$ denotes the space of (equivalence classes of) functions on
$F$ which are square integrable w.r.t. the measure $\mu$.
For the existence of such a map $\tau_F$ see the proof of [39,
thm. 18.6]. By [21, thm. 1, chap. VII] we have that 
$$\mbox{\rm Range}\,\tau_F=H^\alpha(F),\quad
\tau_F\in B(H^s(\RE^n),H^\alpha(F)),\qquad  
\alpha=s-\,\frac{n-d}{2}\ ,
$$
where the Hilbert space $H^\alpha(F)$ is a Besov-like space which
coincides with the usual Sobolev space when $F$ is a regular manifold. In the
case 
$0<\alpha<1$, $H^\alpha(F)$ can
be defined (see [21, \S 1.1, chap. V]) as the set of $f\in L^2(F)$ having finite norm
$$\|f\|^2_{H^\alpha}:=
\|f\|^2_{L^2}+\int_{|x-y|<1}d\mu(x)\,d\mu(y)\,\,
\frac{|f(x)-f(y)|^2}{|x-y|^{d+2\alpha}}\ .
$$
By lemma 2.2, remark 2.12 (taking $\X=H^{\alpha}(F)$ so that
$\tau_F$ is surjective) and theorem 2.1 (hypotheses (h2) being
equivalent 
to $\tau'_F(f)\notin \ld{n}$, $f\not=0$, which is
surely satisfied when $F$ is a null set) one can then immediately
define a family (paramentrized
by the self-adjoint operators on $H^{\alpha}(F)$) of
self-adjoint extensions of $A_{\left|\{\tau_F=0\}\right.}$, where
$A$ is any self-adjoint operator on $L^2(\RE^n)$ with domain $H^s(\RE^n)$.\par  
By considering
$d$-measures one can treat the situation where even more general
sets appear. 
A Borel measure $\mu$ on $\RE^n$ is said to be a
$d$-measure, $d\in (0,n]$, if 
$$
\exists\, c>0\ :\ \forall\, x\in \RE^n,\ \forall\,r\in(0,1],\qquad
\mu(B_r(x))\le cr^d\ .
$$
Then, by [21, lemma 1, chap. VIII], when
$$p=\frac{2d}{n-2s_*},\qquad 0<s_*\le s,\quad n-d<2s_*<n\ ,$$ and
denoting by $L^p(\mu)$ the
space of (equivalence classes of) functions which are $p$-integrable
w.r.t. $\mu$, the linear operator 
$$
\tau_\mu: H^s(\RE^n)\to
L^p(\mu)\qquad\tau_\mu\in B(H^s(\RE^n),L^p(\mu))
$$
$$
f(\tau_\mu\phi):=\int_{\RE^n}d\mu(x)\,f(x)\phi(x)\,,\qquad
f\in L^q(\mu),\ \,\frac{1}{p}+\frac{1}{q}=1\,,
$$
is well defined. Since $\mu_F$, when $F$ is a $d$-set, is a
$d$-measure, 
the previous results tell us that in this case we can
take $n-d=2s_*<2s$ (so that $p=2$) and $\tau_{\mu_F}$ concides with $\tau_F$. 
An interesting example of a $d$-measure is the one
given by the occupation time of Brownian motion: given $\gamma\in
C(\RE_+,\RE^n)$, 
$n\ge 3$, let us define the Radon measure
$$
\mu_\gamma(A):=\int _0^{\infty}dt\,\chi_A(\gamma(t))\ .
$$ 
Then, by estimates on Brownian motion occupation times and by
a Borel-Cantelli argument (see \cite{[CT]}), one has that, for arbitrarily small
positive $\epsilon$ and almost surely with respect
to Wiener measure,
$$
\mu_\gamma(B_r(x))\le cr^{2-\epsilon}\,;
$$
moreover the Hausdorff dimension of the support of $\mu_\gamma$ is equal to two.
\par
Let us now consider the self-adjoint pseudo-differential operator ($s\ge 0$) 
$$\psi(D):H^s(\RE^n)\to L^2(\RE^n)\,,\qquad\psi(D)\phi:=\F^{-1}(\psi\,\F \phi)\ ,
$$
where $\psi$ is a real-valued Borel function such that $$\frac{1}{c}\,(1+|x|^2)^{s/2}\le
1+|\psi(x)|\le c\,(1+ |x|^2)^{s/2}\ .$$
One has
$$
\breve G(z):\ld{n}\to L^p(\mu)\,,\qquad \GB\phi:=\tau_\mu(\K^\psi_z*\phi)\ ,
$$
where 
$$ \K^\psi_z:=\F^{-1}\frac{1}{-\psi+z}\,,\qquad
\K^\psi_z*\phi:=(2\pi)^{-n/2}\,
\F^{-1}\left(\frac{\F\phi}{-\psi+z}\right)
$$
and
$$
G(z):L^q(\mu)\to \ld{n}\,,\qquad G(z)f:=\K^\psi_z*\tau'_\mu(f^*)\ ,
$$
where $\tau'_\mu(f)\in H^{-s}(\RE^n)$ is the signed measure
defined by $$\tau'_\mu(f)(\phi)=\langle f^*,\tau_\mu\phi\rangle\equiv
f(\tau_\mu\phi)$$ and 
$$ 
\K^\psi_z*\tau'_\mu(f):=(2\pi)^{-n/2}\,
\F^{-1}\left(\frac{\F\tau'_\mu(f)}{-\psi+z}\right)\ .
$$
When one uses the representation $\hat\Gamma_\Theta(z)$, by lemma 2.2
one has, if 
$$
\K^\psi_z*\nu(x)=\int_{\RE^n}d\nu(y)\, \K^\psi_z(x-y)\,,\qquad \nu\in H^{-s}(\RE^n)\ ,
$$
and 
$$
\wtilde
\K^\psi_{z_0}:=\frac{1}{2}\,( \K^\psi_{z_0}+
\K^\psi_{z^*_0})\equiv\mbox{\rm Re}(\K^\psi_{z_0})\ ,
$$ 
$$
\hat \Gamma(z): L^q(\mu)\to L^p(\mu)\,,
\qquad \frac{1}{p}+\frac{1}{q}=1\,,
$$
$$
f_1(\hat \Gamma(z)f_2)=
\int_{\RE^{2n}}d\mu(x)\,d\mu(y)\,f_1(x)f^*_2(y)\,
\left(\,\wtilde \K^\psi_{z_0}(x-y)-\K^\psi_z(x-y)\,\right)\ .
$$
By its definition and by the Hahn-Banach theorem we have that $\tau_\mu$ has dense range
when
$$
\left\{f\in L^q(\mu)\ :\ \forall\,\phi\in
H^s(\RE^n)\quad\int_{\RE^n}d\mu(x)\,f(x)\phi(x)=0\right\}=\left\{0\right\}\ .
$$
Therefore $\tau_\mu$ has dense range when the Bessel $s_*$-capacity of
supp($\mu$), $s_*\le s$, is not zero, and this is true ( by
Frostman's lemma, see e.g. [29, thm. 7.1]) when
$$
n-d(\mu)<2s_*\le n,\qquad 0<s_*\le s\ ,
$$ 
where $d(\mu)$ denotes the Hausdorff
dimension of supp$(\mu)$. \par 
Let us note that, since $p\in
[2,\infty)$ in (19), when $\mu$ is a finite measure we can view $\tau_\mu$ as
a map into the Hilbert space $L^2(\mu)$. In this  case we can then try
to apply lemma 2.4 in order to find other maps $\Gamma(z)$ which satisfy
(5) and (7.1). Supposing that $\psi(x)=\psi(-x)$, so that
$\K^\psi_z(x-y)=\K^\psi_z(y-x)$, and  
$$
\int_{\RE^{2n}}d\mu(x)\,d\mu(y)\,|x-y|^2\,|\K^\psi_z(x-y)|\ <+\infty
$$
we have, for any $f_1,f_2\in C_0^1(\RE^n)$, and proceeding similarly to
example 3.5 (case $n=4$),
\begin{eqnarray*}
&\quad&\langle f_1,\hat \Gamma(z)f_2\rangle\\
&=&\int_{\RE^{2n}}d\mu(x)\,d\mu(y)\,f^*_1(x)f_2(y)\,
\left(\,\wtilde \K^\psi_{z_0}(x-y)-\K^\psi_z(x-y)\,\right)\cr
&=&\frac{1}{2}\int_{\RE^{2n}}d\mu(x)\,d\mu(y)\,(\,f^*_1(x)-f^*_1(y)\,)\,
(\,f_2(x)-f_2(y)\,)\,\K^\psi_{z}(x-y)\cr
&\quad&+\int_{\RE^{2n}}d\mu(x)\,d\mu(y)\,f^*_1(x)f_2(x)\,
\left(\,\wtilde \K^\psi_{z_0}(x-y)-\K^\psi_z(x-y)\,\right)\cr
&\quad&-\frac{1}{2}\int_{\RE^{2n}}d\mu(x)\,d\mu(y)\,(\,f^*_1(x)-f^*_1(y)\,)\,
(\,f_2(x)-f_2(y)\,)\,\wtilde\K^\psi_{z_0}(x-y)
\ .
\end{eqnarray*}
Therefore, with the last term being $z$-independent, the sequilinear form
$$
\wtilde\E(z):C^1_0(\RE^n)\times C^1_0(\RE^n)\to\C
$$
\begin{eqnarray*}
&\quad&\wtilde\E(z)(f_1,f_2)\\
&=&\frac{1}{2}\int_{\RE^{2n}}d\mu(x)\,d\mu(y)\,(\,f^*_1(x)-f^*_1(y)\,)\,
(\,f_2(x)-f_2(y)\,)\,\K^\psi_{z}(x-y)\cr
&\quad&+\int_{\RE^{2n}}d\mu(x)\,d\mu(y)\,f^*_1(x)f_2(x)\,
\left(\,\wtilde \K^\psi_{z_0}(x-y)-\K^\psi_z(x-y)\,\right)
\end{eqnarray*}
satisfies (13) and (14). In the case $\wtilde\K^\psi_{z_0}\ge 0$ one
has $\mbox{\rm Re}\,(\wtilde\E(z_0)(f,f))\ge 0$ 
and so (15) is satisfied with $M=0$. Moreover $\wtilde\E(z_0)$ is
readily checked to 
be closable. So, by lemma 2.4, proposition 2.1 and
theorem 2.1, for any strictly positive (i.e. $\gamma(\Theta)>0$) 
self-adjoint operator on $L^2(\mu)$ one obtains a family of
self-adjoint 
extensions $\psi(D)_\Theta^\mu$
with resolvent
$$
(-\psi(D)_\Theta^\mu+z)^{-1}\phi
=(-\psi(D)+z)^{-1}\phi+\K^\psi_z*\tau_\mu'(\,\Theta+
\Gamma(z)^{-1}\cdot\tau_\mu(\K^\psi_z*\phi)\,)\ ,
$$
where $\Gamma(z)$ is the operator corresponding to the closure of
$\wtilde\E(z)$. Such a family, in the particular case 
$\psi(D)=\Delta$, is the same obtained, by an
approximation method, in \cite{[AFHKL]} (also see \cite{[C]}) and 
generalizes, although with a
different $\Gamma(z)$, the situation
discussed in example 3.5. In this regard suppose that the 
$1$-set $C$ is the range of a Lipschitz path
$\gamma:I\subseteq\RE\to C\subset\RE^n$, $n=3$ or $n=4$, $|\dot\gamma|=1$ a.e., 
(so that $\tau_{\mu_C}$ has dense range
in $L^2(C)\simeq L^2(I)$). Under hypothesis (18) one can
again consider, when $n=3$ the operator $\wtilde \Gamma(z)$ appearing in (19) 
and, when $n=4$ the sesquilinear form $\wtilde\E(z)$ appearing in (20), the only
difference being that now the domain of definition of such objects is 
$C^1_0(I\backslash I_*)$, with 
\begin{eqnarray*}
I_*&:=&\left\{t\in I\, :\, \mbox{\rm $\gamma$ is not differentiable at $t$}\right\}\\
&\cup&\left\{t\in I\, :\, \exists \,s\not=t\quad \mbox{\rm s.t.}\quad
\gamma(t)=\gamma(s)\right\}
\end{eqnarray*}
(of course, in order $C^1_0(I\backslash I_*)$ to be still a dense set,
one has to suppose that the closure of $I_*$ is a null set).
However in the case $n=4$ the problem of the semi-boundedness of $\wtilde\E(z)$
arises: indeed one can show (see \cite{[Sh]}) that $\wtilde\E(z)$ is unbounded
from below in the case where $\gamma$ has angle points. This phenomenon is
similar to the one related to unboundedness from below of
Schr\"odinger operators describing $n(>2)$ point interacting particles 
(see \cite{[MF]}, \cite{[DFT]} and references therein). 
\end{example}
\begin{example} {\it Singular perturbations of the d'Alembertian
supported by time-like straight lines.} We take $\calH=L^2(\RE^4)$,
$$A=\square:=-\Delta_{(1)}\otimes\uno+\uno\otimes\Delta_{(3)}\,,$$
$\Delta_{(d)}$ being the Laplacian in $d$ dimensions, 
and ($h\in\RE$, $k\in\RE^3$ denoting the variables dual to $t\in\RE$, $x\in\RE^3$)
$$D(\square)=\left\{\Phi\in L^2(\RE^4)\ :\
(h^2-|k|^2)\F\Phi(h,k)\in L^2(\RE^4)\right\}\ .$$
Let $\ell(s)=y+ws$,
$y,w\in\RE^4$, 
be a time-like straight line, i.e. 
$$
w=(\gamma_v,\gamma_v\,v),\quad v\in\RE^3,\quad |v|<1,\quad \gamma_v:=
\left(1-|v|^2\right)^{-1/2}\ .
$$
Consider now the unique surjective linear operator 
$$
\tau_0: D(\square)\to H^{-{1}/{2}}(\RE)\,,\qquad\tau_0\in B(D(\square),
H^{-{1}/{2}}(\RE))
$$
such that
$$
\forall\Phi\in C^\infty_0(\RE^4),\qquad\tau_0\Phi(s):=\Phi(s,0)\ .
$$
For the existence of such a $\tau_0$ see the next example. \par 
Let then $\Pi_{y,v}$ 
be the unitary operator which compose any
function in $L^2(\RE^4)$ with the Lorentz boost corresponding to $v$ and then with
the translation by $y$, so that $\Pi_{y,v}\in B(D(\square),D(\square))$. Defining 
$$
\tau_{y,v}:=\tau_0\cdot\Pi_{y,v}:D(\square)\to
H^{-{1}/{2}}(\RE)\,,\qquad\tau_{y,v}\in B(D(\square),
H^{-{1}/{2}}(\RE))
$$
one has 
$$
\forall\Phi\in C^\infty_0(\RE^4),\qquad\tau_{y,v}\Phi(s):=\Phi(\ell(s))\ .
$$
We begin studying the self-adjoint extensions given by $\tau_0$. By
Fourier transform (here and below $z\in\C\backslash\RE$) one obviously has
$$
\F\cdot(-\square+z)^{-1}\Phi(h,k)=\frac{\F\Phi(h,k)}{-h^2+|k|^2+z}\ .
$$
So, since, as $R\uparrow\infty$,
$$
\frac{1}{\sqrt {2\pi}}\int_0^\infty dh\,\frac{1}{\left|-h^2+R^2+z\right|^2}\
\sim\frac{c}{4R}\ ,
$$ 
by the H\"older inequality and the Riemann-Lebesgue lemma there follows that
$$
\forall\,\Phi=\phi\otimes\varphi\in L^2(\RE)\otimes H^s(\RE^3),\ s>1,\qquad 
(-\square+z)^{-1}\Phi\in C_b^0(\RE^4)
$$
and, by Fubini theorem,
\begin{eqnarray*}
&\quad&[(-\square+z)^{-1}\Phi]\,(t,x)\cr
&=&{1\over{(2\pi)^2}}\,\int_{\RE^4}dh\,dk\,e^{iht}e^{ik\cdot
x}{{\F\Phi (h,k)}\over{-h^2+|k|^2+z}}\cr
&=&{1\over{(2\pi)^{1/2}}}\,\int_{\RE}dh\,e^{iht}\F\phi(h)
\left(\left(-\Delta-h^2+z\right)^{-1}\varphi\right)(x)\cr
&=&{1\over{(2\pi)^{1/2}}}\,\int_{\RE}dh\,e^{iht}\F\phi(h)
\int_{\RE^3}dy\,\vp(y)\,\frac{e^{-\sqrt{-h^2+z}\,|x-y|}}{4\pi\,|x-y|}\cr
&=&
\int_{\RE^3}dy\,\frac{\vp(y)}{4\pi\,|x-y|}\,
\left(e^{-|x-y|\,\sqrt{\Delta_{(1)}+z}}\,\phi\right)(t)
\ .
\end{eqnarray*}
Here $\mbox{\rm Re}\sqrt{-h^2+z}>0$; this choice will be always
assumed in the following without further specification. The above
calculation then gives 
$$
\breve G(z):L^2(\RE)\otimes H^{s}(\RE^3)\to C_b^0({\RE})\, ,\quad s>1,
$$
$$
\left(\breve G(z)\phi\otimes\varphi\right)\,(t):=
\int_{\RE^3}dy\,\frac{\varphi(y)}{4\pi\,|y|}\,
\left[e^{-|y|\,\sqrt{\Delta_{(1)}+z}}\phi\right](t)
$$
and
$$
G(z):H^{1/2}(\RE)\to L^2(\RE^4)$$
$$
(G(z)\phi)\,(t,x):=\frac{1}{4\pi\,|x|}\,
\left[e^{-|x|\,\sqrt{\Delta_{(1)}+z}}\phi^*\right](t)\ .
$$
Let us note that $\breve G(z)$ extends to a continuous linear operator
from $L^2(\RE^4)$ to $H^{-1/2}(\RE)$ since
\begin{eqnarray*}
&\quad&\|\breve G(z)\phi\otimes\varphi\|^2_{H^{-1/2}}\cr
&\le&\,
\int_\RE dt\left(\,\int_{\RE^3}dy\,\frac{|\varphi(y)|}{4\pi\,|y|}\,
\left|\left[(-\Delta_{(1)}+1)^{-1/4}\cdot e^{-|y|\,\sqrt{\Delta_{(1)}+z}}
\phi\right](t)\right|\,\right)^2\cr
&\le&\,\|\varphi\|_{L^2}^2
\,\left\langle\left(-\Delta_{(1)}+1\right)^{-1/2}\cdot
\int_0^\infty dR\,
e^{-2R\,{\mathrm {Re}}\left(\sqrt{\Delta_{(1)}+z}\,\right)}\phi,\phi\right\rangle\cr
&=&\,\frac{1}{2}\,\|\varphi\|_{L^2}^2
\,\left\langle\left(-\Delta_{(1)}+1\right)^{-1/2}\cdot
\left(\mbox{\rm
Re}\left(\sqrt{\Delta_{(1)}+z}\,\right)\right)^{-1}\phi,\phi
\right\rangle\cr
&\le&\,\frac{1}{2}\,\|\phi\otimes\varphi\|_{L^2}^2
\,\left\|\left(\mbox{\rm
Re}\left(\sqrt{\Delta_{(1)}+z}\,\right)\right)^{-1/2}\right\|_{L^2,H^{-1/2}}\ .
\end{eqnarray*}
Similarly $G(z)$ is a continuous linear operator from $H^{1/2}(\RE)$ to
$L^2(\RE^4)$ 
since 
\begin{eqnarray*}
&\quad&\|G(z)\phi\|^2_{L^2}\cr
&=&\int_0^\infty dR\,\left\| 
e^{-R\,\sqrt{\Delta_{(1)}+z}}
\phi\right\|^2_{L^2}\cr
&=&\,\left\langle\int_0^\infty dR\,
e^{-2R\,{\mathrm
{Re}}\left(\sqrt{\Delta_{(1)}+z}\,\right)}\phi,\phi
\right\rangle\cr
&=&\,\frac{1}{2}
\,\left\langle
\left(\mbox{\rm
Re}\left(\sqrt{\Delta_{(1)}+z}\,\right)\right)^{-1}\phi,\phi
\right\rangle\cr
&\le&\,\frac{1}{2}\,\|\phi\|_{H^{1/2}}^2
\,\left\|\left(\mbox{\rm
Re}\left(\sqrt{\Delta_{(1)}+z}\,\right)\right)^{-1/2}\right\|_{H^{1/2},L^2}\ .
\end{eqnarray*}
We now look for the map $\Gamma(z)$. Since
\begin{eqnarray*}
&\quad&(z-w)\,\F\cdot\breve G(w)\cdot
G(z)\phi(h)\cr
&=&
\frac{z-w}{(2\pi)^3}\,\F\phi^*(h)\int_{\RE^3} dk\,
\frac{1}{(-h^2+|k|^2+w)(-h^2+|k|^2+z)}\cr
&=&\frac{z-w}{2\pi^2}\,\F\phi^*(h)\int_0^\infty dr\,
\frac{r^2}{(-h^2+r^2+w)(-h^2+r^2+z)}\cr
&=&\frac{1}{4\pi}\,\left(\sqrt{-h^2+z}-\sqrt{-h^2+w}\,\right)\F\phi^*(h)
\end{eqnarray*}
one defines
$$
\Gamma(z)\ :H^{1/2}(\RE)\to H^{-1/2}(\RE)
$$
$$
\Gamma(z)\phi:=\frac{1}{4\pi}\,\sqrt{\Delta_{(1)}+z}\ \phi^*\ .
$$
Of course we can view $\Gamma(z)$ as a (unbounded) closed and densely
defined linear operator on the Hilbert space
$H^{-1/2}(\RE)$; evidently $\Gamma(z)^*=\Gamma(z)$. 
Therefore, by proposition 2.1, $\Gamma_\Theta(z)$
satisfies (h1) (with $Z_\Theta=\rho(\square)$) for any self-adjoint 
operator $\Theta$ on $H^{-1/2}(\RE)$ which is $\Gamma(z)$-bounded. 
It is immediate, by Fourier transform, to check the
validity of (h2). Therefore the trace $\tau_0$ gives rise to the family of self-adjoint
extensions $\square^0_\Theta$ with resolvent
$$
(-\square^0_\Theta+z)^{-1}=(-\square+z)^{-1}
+G(z)\cdot\left(\,\Theta+\frac{1}{4\pi}\,\sqrt{\Delta_{(1)}+z}\,\right)^{-1}
\cdot\breve G(z)
$$ 
(here, since they annihilates between themselves, we did not put the 
complex conjugations appearing in both the definitions of $G(z)$ and 
$\Gamma_\Theta(z)$).
By our definition of $\tau_{y,v}$ we have, since $\Pi_{y,v}$ commutes
with $\square$, 
$$
\breve G_{y,v}(z):=\tau_{y,v}\cdot R(z)=\breve G(z)\cdot\Pi_{y,v}
$$
and
$$
C_{L^2}^{-1}\cdot\breve G_{y,v}(z^*)'=\Pi^*_{y,v}\cdot
G(z)\ .
$$
This immediately implies that one can use the same $\Gamma_\Theta(z)$
as before and so the trace $\tau_{y,v}$ gives rise to the family of self-adjoint
extensions $\square^{y,v}_\Theta$ with resolvent
\begin{eqnarray*}
&\quad&(-\square^{y,v}_\Theta+z)^{-1}\\
&=&(-\square+z)^{-1}
+\Pi^*_{y,v}\cdot
G(z)\cdot\left(\,\Theta+\frac{1}{4\pi}\,\sqrt{\Delta_{(1)}+z}\,\right)^{-1}
\cdot\breve G(z)\cdot\Pi_{y,v}\ .
\end{eqnarray*}
Moreover the following kind of Poincar\'e-invariance holds:
$$
D(\square^{y,v}_\Theta)=\Pi^*_{y,v}(D(\square^0_{\Theta}))\ \quad
\mbox{\rm and}\quad\ 
\square^{y,v}_\Theta=\Pi^*_{y,v}\cdot\square^0_{\Theta}\cdot
\Pi_{y,v}\ .
$$
Let us remark that, even if the operator
$\Gamma_\Theta(z)$ appearing in the resolvent above coincides with the one
used in the case $v=0$, it is 
applied to functions which depend on different variables: 
when $v=0$ it acts
on functions of the relative time whereas it acts on functions of the
proper time when $v\not= 0$. Therefore if in the case $v\not= 0$ one
uses relative time, then
$\Gamma_\Theta(z)$ becomes a velocity-dependent operator.  
\end{example}
\begin{example}{\it Singular perturbation given by traces on
Malgrance spaces.} Given any continuous functions 
$\varphi>0$ on $\RE^n$, $\varphi\in\M$ will mean that there
exists a polynomial $P$ such that
$$
\forall\,x\in\RE^n,\qquad \frac{1}{|P(x)|}\le \varphi(x)\le |P(x)|\ .
$$
Then we define the Hilbert space
$H_\varphi(\RE^n)$, $\varphi\in\M$, as the set of tempered distribution $f$ such that
$\F f$ is a functions and
$$
\|f\|^2_\varphi:=\int_{\RE^n}|\varphi(k)\F f(k)|^2\, dk<+\infty\,.
$$
Such a class of function spaces were introduced by Malgrange
in \cite{[M]}.\par
In connection with the previous examples note that 
$$\varphi(x)=(1+|x|^2)^{s/2},\ s\in\RE\quad\Longrightarrow\quad
H_\varphi(\RE^n)=H^s(\RE^n)$$ 
and $$
\varphi(t,x)=(1+(-t^2+|x|^2)^2)^{1/2},\ t\in\RE,\ x\in\RE^3
\quad\Longrightarrow\quad H_\varphi(\RE^4)=D(\square)\ .
$$
We list now some properties of the spaces $H_\varphi(\RE^n)$ following 
[28, \S 1], [19, \S II.2.2] and \cite{[VP]}. Let us remark here that the definition of
$H_\varphi(\RE^n)$ given in \cite{[H]} and \cite{[VP]} is different: it corresponds to
the case in which $\varphi$ belongs to the narrower class $\K$
defined by
$$
\varphi\in\K\ \iff\ \exists\, c,\,N>0\ :\ \forall\,x,\,y\in\RE^n,\qquad 
\varphi(x+y)\le (\,1+c\,|x|\,)^N\varphi(y)\ .
$$
The choice $\varphi\in\K$ ensures that $H_\varphi(\RE^n)$ is a
module over $C^\infty_0(\RE^n)$. However
the results we will quote from [19, \S II.2.2] and \cite{[VP]} hold true also
for the more general case in which $\varphi\in\M$ (see [41, remark 2.3]).\par
The dual space of $H_\vp$ can be
explicitly characterized (see [28, \S 1.1], [19,
thm. 2.2.9], [41, \S 2.1]) as
$$
H_\vp(\RE^n)'\simeq H_{1/\vp}(\RE^n)\ .
$$
As regards the relation between different spaces, by [19, thm. 2.2.2] one has
$$
\vp_1\le c\,\vp_2\quad\Longrightarrow
\quad H_{\varphi_1}(\RE^n)\subseteq
H_{\varphi_2}(\RE^n)\ ,
$$
the embedding being continuous. Therefore, for any $\vp\in\M$ such that 
$\vp\ge c>0$, one has
$H_\vp(\RE^n)\subseteq L^2(\RE^n)$; $H_\vp(\RE^n)$ is then dense in $L^2(\RE^n)$
since $C^\infty_0(\RE^n)$ is dense in $H_\vp(\RE^n)$ (see [28, \S 1.1],
[19, thm. 2.2.1], [41, \S 2.1]). The regularity of elements in
$H_\vp(\RE^n)$ is given by [19, thm. 2.2.7]:
$$
(1+|x|)^k/\vp(x)\in L^2(\RE^n)\quad\Longrightarrow\quad
H_\vp(\RE^n)\subset C^k(\RE^n)\ ,
$$ 
the embedding being continuous. \par
Let us now come to the trace operator on $H_\vp(\RE^n)$ (see [19, thm. 2.2.8],
[41, \S 6]). We write $\RE^n=\RE^d\oplus\RE^{n-d}$, $1\le
d\le n-1$, $x=(\tilde x,\hat x)$, $\tilde x\in\RE^d$, $\hat
x\in\RE^{n-d}$. Suppose that 
$$
\left(\int_{\RE^{n-d}} \frac{1}{\vp^2(0,\hat x)}\,d\hat x\right)^{-1/2}\,<+\infty\ . 
$$  
Then there exists an unique surjective linear operator $\tau_{(d)}$
$$
\tau_{(d)}:H_\vp(\RE^n)\to H_{\tilde\vp}(\RE^d)\,,\qquad \tau_{(d)}\in
B(H_\vp(\RE^n),H_{\tilde\vp}(\RE^d))\ ,
$$
$$
\tilde\vp(\tilde x):=\left(\int_{\RE^{n-d}} 
\frac{1}{\vp^2(\tilde x,\hat x)}\,d\hat x\right)^{-1/2}\ ,
$$
such that 
$$
\forall\,\phi\in C_0^\infty(\RE^n)\,,\qquad\tau_{(d)}(\phi)(\tilde
x)=\phi(\tilde x,0)\ .
$$
The reader can check that the case
$\varphi(t,x)=(1+(-t^2+|x|^2)^2)^{1/2}$, $d=1$, reproduces the trace $\tau_0$
given in the previous example. \par
The trace $\tau_{(d)}$ can be generalized to cover the case of
non-linear 
subsets in the
following way: let $\mu\in H'_\phi(\RE^n)$, $\phi\in\K$ 
(for example, $\mu$ could
be the 
Hausdorff measure of some subset of $\RE^n$ but more general
distributions are allowed), for which there exists $\tilde\phi\in\K$
such that
$$
\int_{\RE^n}\frac{\phi^2(x-y)}{\varphi^2(x)\,\tilde\phi^2(y)}\,dy<c<+\infty\ .
$$ 
Then, by [41, \S 2.4], 
$$\forall f\in H_{\tilde\phi}(\RE^n),
\qquad f\mu\in H'_\vp(\RE^n)\
,$$
where 
$$
f\mu:=(2\pi)^{-n/2}\F^{-1}(\F f*\F \mu)\ .
$$
So we can define
$$
\tilde\tau_\mu:H_{\tilde\phi}(\RE^n)\to
H_\vp'(\RE^n)\qquad\tilde\tau_\mu(f):=f\mu\ ,
$$
and then we have a trace generalizing $\tau_{(d)}$ by 
$$
\tau_\mu:H_\vp(\RE^n)\to H_{\tilde\phi}'(\RE^n)\qquad
\tau_\mu:=\tilde\tau_\mu'\cdot J_{H_\vp}\ .
$$
Let us now consider the self-adjoint pseudo-differential operator
(here \p$\vp\ge c>0$) 
$$\psi(D):H_\vp(\RE^n)\to L^2(\RE^n)\,,\qquad
\psi(D)\Phi:=\F^{-1}(\psi\,\F \Phi)\ ,
$$
where $\psi$ is a real-valued Borel function such that $$\frac{1}{c}\,\vp(x)\le
1+|\psi(x)|\le c\,\vp(x)\ .$$
By Fourier transform one has, if $\tau_\mu$ is defined as above,
$$
G(z)f=G^\mu(z)f:=\frac{1}{(2\pi)^{n/2}}\,\F^{-1}\left(\frac{\F f^**\F\mu}{-\psi+z^*}\right)\ .
$$
Therefore (h2) is equivalent to $\F f*\F\mu\notin \ld{n}$,
i.e. $f\mu\notin\ld{n}$. This condition is surely satisfied when the
support of $\mu$ is a set of zero Lebesgue measure. \par
By lemma 2.2 we have then, for any $f_1,f_2\in H_{\tilde\phi}(\RE^n)$, 
\begin{equation}
\hat\Gamma(z)f_1(f_2)=f_2\mu((\wtilde G^\mu-G^\mu(z))f_1)\ , 
\end{equation}
where 
$$
\wtilde G^\mu:=\frac{G^\mu(z_0)+G^\mu(z_0^*)}{2}\,,\qquad z_0\in\rho(\psi(D))\ .
$$
In the case  
$$
\vp\in\K,\qquad\int_{\RE^n}\frac{\phi^2(x-y)}{\tilde\phi^2(x)\,\vp^2(y)}\,dy<c<+\infty\ ,
$$
by [41, \S 2.4] as above, we have  
$$\forall \Phi\in H_{\vp}(\RE^n),
\qquad \Phi\mu\in H'_{\tilde\phi}(\RE^n)\ ,$$
and so, by (22), $\hat\Gamma(z)=\hat\Gamma^\mu(z)$, where
$$
\hat\Gamma^\mu(z):H_{\tilde\phi}''(\RE^n)\simeq H_{\tilde\phi}(\RE^n)
\to H_{\tilde\phi}'(\RE^n)\,,\qquad
\hat\Gamma^\mu(z)f:=((\wtilde G^\mu-G^\mu(z))f)\mu\ .
$$ 
In the case where $\tau_\mu$ is surjective, by remark 2.12,
proposition 2.1 and theorem 2.1, $\tau_\mu$ gives rise to the family of
self-adjoint operators $\psi(D)^\mu_\Theta$ with resolvents
$$
(-\psi(D)^\mu_\Theta+z)^{-1}=(-\psi(D)+z)^{-1}+G^{\mu}(z)\cdot
(\Theta+\hat\Gamma^\mu(z))^{-1}\cdot\breve G^\mu(z)\ ,
$$
where $\Theta$ is any operator from  
$H_{\tilde\phi}(\RE^n)$ to $H_{\tilde\phi}'(\RE^n)$ such that $\Theta=\Theta'$.
\end{example}


\begin{acknowledgment}
We thank Sergio Albeverio, Gianfausto Dell'Antonio, Diego Noja and
Sandro Teta for stimulating discussions.
\end{acknowledgment}


\end{article}
\end{document}